\documentclass[reqno,12pt,a4paper]{amsart}

\usepackage{mathrsfs}
\usepackage{amssymb}
\usepackage{amsfonts}
\usepackage{latexsym}
\usepackage{amsthm}
\usepackage{graphicx}

\def\lb{\label}

\newcommand{\er}[1]{\textrm{(\ref{#1})}}

\begin{document}


\renewcommand{\theequation}{\arabic{section}.\arabic{equation}}
\theoremstyle{plain}
\newtheorem{theorem}{\bf Theorem}[section]
\newtheorem{lemma}[theorem]{\bf Lemma}
\newtheorem{corollary}[theorem]{\bf Corollary}
\newtheorem{proposition}[theorem]{\bf Proposition}
\newtheorem{definition}[theorem]{\bf Definition}

\newtheorem{remark}[theorem]{\bf Remark}

\def\a{\alpha}  \def\cA{{\mathcal A}}     \def\bA{{\bf A}}  \def\mA{{\mathscr A}}
\def\b{\beta}   \def\cB{{\mathcal B}}     \def\bB{{\bf B}}  \def\mB{{\mathscr B}}
\def\g{\gamma}  \def\cC{{\mathcal C}}     \def\bC{{\bf C}}  \def\mC{{\mathscr C}}
\def\G{\Gamma}  \def\cD{{\mathcal D}}     \def\bD{{\bf D}}  \def\mD{{\mathscr D}}
\def\d{\delta}  \def\cE{{\mathcal E}}     \def\bE{{\bf E}}  \def\mE{{\mathscr E}}
\def\D{\Delta}  \def\cF{{\mathcal F}}     \def\bF{{\bf F}}  \def\mF{{\mathscr F}}
\def\c{\chi}    \def\cG{{\mathcal G}}     \def\bG{{\bf G}}  \def\mG{{\mathscr G}}
\def\z{\zeta}   \def\cH{{\mathcal H}}     \def\bH{{\bf H}}  \def\mH{{\mathscr H}}
\def\e{\eta}    \def\cI{{\mathcal I}}     \def\bI{{\bf I}}  \def\mI{{\mathscr I}}
\def\p{\psi}    \def\cJ{{\mathcal J}}     \def\bJ{{\bf J}}  \def\mJ{{\mathscr J}}
\def\vT{\Theta} \def\cK{{\mathcal K}}     \def\bK{{\bf K}}  \def\mK{{\mathscr K}}
\def\k{\kappa}  \def\cL{{\mathcal L}}     \def\bL{{\bf L}}  \def\mL{{\mathscr L}}
\def\l{\lambda} \def\cM{{\mathcal M}}     \def\bM{{\bf M}}  \def\mM{{\mathscr M}}
\def\L{\Lambda} \def\cN{{\mathcal N}}     \def\bN{{\bf N}}  \def\mN{{\mathscr N}}
\def\m{\mu}     \def\cO{{\mathcal O}}     \def\bO{{\bf O}}  \def\mO{{\mathscr O}}
\def\n{\nu}     \def\cP{{\mathcal P}}     \def\bP{{\bf P}}  \def\mP{{\mathscr P}}
\def\r{\rho}    \def\cQ{{\mathcal Q}}     \def\bQ{{\bf Q}}  \def\mQ{{\mathscr Q}}
\def\s{\sigma}  \def\cR{{\mathcal R}}     \def\bR{{\bf R}}  \def\mR{{\mathscr R}}
\def\S{\Sigma}  \def\cS{{\mathcal S}}     \def\bS{{\bf S}}  \def\mS{{\mathscr S}}
\def\t{\tau}    \def\cT{{\mathcal T}}     \def\bT{{\bf T}}  \def\mT{{\mathscr T}}
\def\f{\phi}    \def\cU{{\mathcal U}}     \def\bU{{\bf U}}  \def\mU{{\mathscr U}}
\def\F{\Phi}    \def\cV{{\mathcal V}}     \def\bV{{\bf V}}  \def\mV{{\mathscr V}}
\def\P{\Psi}    \def\cW{{\mathcal W}}     \def\bW{{\bf W}}  \def\mW{{\mathscr W}}
\def\o{\omega}  \def\cX{{\mathcal X}}     \def\bX{{\bf X}}  \def\mX{{\mathscr X}}
\def\x{\xi}     \def\cY{{\mathcal Y}}     \def\bY{{\bf Y}}  \def\mY{{\mathscr Y}}
\def\X{\Xi}     \def\cZ{{\mathcal Z}}     \def\bZ{{\bf Z}}  \def\mZ{{\mathscr Z}}
\def\be{{\bf e}}
\def\bv{{\bf v}} \def\bu{{\bf u}}
\def\Om{\Omega}
\def\bbD{\pmb \Delta}
\def\mm{\mathrm m}
\def\mn{\mathrm n}

\newcommand{\mc}{\mathscr {c}}

\newcommand{\gA}{\mathfrak{A}}          \newcommand{\ga}{\mathfrak{a}}
\newcommand{\gB}{\mathfrak{B}}          \newcommand{\gb}{\mathfrak{b}}
\newcommand{\gC}{\mathfrak{C}}          \newcommand{\gc}{\mathfrak{c}}
\newcommand{\gD}{\mathfrak{D}}          \newcommand{\gd}{\mathfrak{d}}
\newcommand{\gE}{\mathfrak{E}}
\newcommand{\gF}{\mathfrak{F}}           \newcommand{\gf}{\mathfrak{f}}
\newcommand{\gG}{\mathfrak{G}}           
\newcommand{\gH}{\mathfrak{H}}           \newcommand{\gh}{\mathfrak{h}}
\newcommand{\gI}{\mathfrak{I}}           \newcommand{\gi}{\mathfrak{i}}
\newcommand{\gJ}{\mathfrak{J}}           \newcommand{\gj}{\mathfrak{j}}
\newcommand{\gK}{\mathfrak{K}}            \newcommand{\gk}{\mathfrak{k}}
\newcommand{\gL}{\mathfrak{L}}            \newcommand{\gl}{\mathfrak{l}}
\newcommand{\gM}{\mathfrak{M}}            \newcommand{\gm}{\mathfrak{m}}
\newcommand{\gN}{\mathfrak{N}}            \newcommand{\gn}{\mathfrak{n}}
\newcommand{\gO}{\mathfrak{O}}
\newcommand{\gP}{\mathfrak{P}}             \newcommand{\gp}{\mathfrak{p}}
\newcommand{\gQ}{\mathfrak{Q}}             \newcommand{\gq}{\mathfrak{q}}
\newcommand{\gR}{\mathfrak{R}}             \newcommand{\gr}{\mathfrak{r}}
\newcommand{\gS}{\mathfrak{S}}              \newcommand{\gs}{\mathfrak{s}}
\newcommand{\gT}{\mathfrak{T}}             \newcommand{\gt}{\mathfrak{t}}
\newcommand{\gU}{\mathfrak{U}}             \newcommand{\gu}{\mathfrak{u}}
\newcommand{\gV}{\mathfrak{V}}             \newcommand{\gv}{\mathfrak{v}}
\newcommand{\gW}{\mathfrak{W}}             \newcommand{\gw}{\mathfrak{w}}
\newcommand{\gX}{\mathfrak{X}}               \newcommand{\gx}{\mathfrak{x}}
\newcommand{\gY}{\mathfrak{Y}}              \newcommand{\gy}{\mathfrak{y}}
\newcommand{\gZ}{\mathfrak{Z}}             \newcommand{\gz}{\mathfrak{z}}

\def\ve{\varepsilon}   \def\vt{\vartheta}    \def\vp{\varphi}    \def\vk{\varkappa}

\def\A{{\mathbb A}} \def\B{{\mathbb B}} \def\C{{\mathbb C}}
\def\dD{{\mathbb D}} \def\E{{\mathbb E}} \def\dF{{\mathbb F}} \def\dG{{\mathbb G}} \def\H{{\mathbb H}}\def\I{{\mathbb I}} \def\J{{\mathbb J}} \def\K{{\mathbb K}} \def\dL{{\mathbb L}}\def\M{{\mathbb M}} \def\N{{\mathbb N}} \def\O{{\mathbb O}} \def\dP{{\mathbb P}} \def\R{{\mathbb R}}\def\S{{\mathbb S}} \def\T{{\mathbb T}} \def\U{{\mathbb U}} \def\V{{\mathbb V}}\def\W{{\mathbb W}} \def\X{{\mathbb X}} \def\Y{{\mathbb Y}} \def\Z{{\mathbb Z}}


\def\la{\leftarrow}              \def\ra{\rightarrow}            \def\Ra{\Rightarrow}
\def\ua{\uparrow}                \def\da{\downarrow}
\def\lra{\leftrightarrow}        \def\Lra{\Leftrightarrow}


\def\lt{\biggl}                  \def\rt{\biggr}
\def\ol{\overline}               \def\wt{\widetilde}
\def\ul{\underline}
\def\no{\noindent}


\let\ge\geqslant                 \let\le\leqslant
\def\lan{\langle}                \def\ran{\rangle}
\def\/{\over}                    \def\iy{\infty}
\def\sm{\setminus}               \def\es{\emptyset}
\def\ss{\subset}                 \def\ts{\times}
\def\pa{\partial}                \def\os{\oplus}
\def\om{\ominus}                 \def\ev{\equiv}
\def\iint{\int\!\!\!\int}        \def\iintt{\mathop{\int\!\!\int\!\!\dots\!\!\int}\limits}
\def\el2{\ell^{\,2}}             \def\1{1\!\!1}
\def\sh{\sharp}
\def\wh{\widehat}
\def\bs{\backslash}
\def\intl{\int\limits}

\def\na{\mathop{\mathrm{\nabla}}\nolimits}
\def\sh{\mathop{\mathrm{sh}}\nolimits}
\def\ch{\mathop{\mathrm{ch}}\nolimits}
\def\where{\mathop{\mathrm{where}}\nolimits}
\def\all{\mathop{\mathrm{all}}\nolimits}
\def\as{\mathop{\mathrm{as}}\nolimits}
\def\Area{\mathop{\mathrm{Area}}\nolimits}
\def\arg{\mathop{\mathrm{arg}}\nolimits}
\def\const{\mathop{\mathrm{const}}\nolimits}
\def\det{\mathop{\mathrm{det}}\nolimits}
\def\diag{\mathop{\mathrm{diag}}\nolimits}
\def\diam{\mathop{\mathrm{diam}}\nolimits}
\def\dim{\mathop{\mathrm{dim}}\nolimits}
\def\dist{\mathop{\mathrm{dist}}\nolimits}
\def\Im{\mathop{\mathrm{Im}}\nolimits}
\def\Iso{\mathop{\mathrm{Iso}}\nolimits}
\def\Ker{\mathop{\mathrm{Ker}}\nolimits}
\def\Lip{\mathop{\mathrm{Lip}}\nolimits}
\def\rank{\mathop{\mathrm{rank}}\limits}
\def\Ran{\mathop{\mathrm{Ran}}\nolimits}
\def\Re{\mathop{\mathrm{Re}}\nolimits}
\def\Res{\mathop{\mathrm{Res}}\nolimits}
\def\res{\mathop{\mathrm{res}}\limits}
\def\sign{\mathop{\mathrm{sign}}\nolimits}
\def\span{\mathop{\mathrm{span}}\nolimits}
\def\supp{\mathop{\mathrm{supp}}\nolimits}
\def\Tr{\mathop{\mathrm{Tr}}\nolimits}
\def\BBox{\hspace{1mm}\vrule height6pt width5.5pt depth0pt \hspace{6pt}}


\newcommand\nh[2]{\widehat{#1}\vphantom{#1}^{(#2)}}
\def\dia{\diamond}

\def\Oplus{\bigoplus\nolimits}



\def\qqq{\qquad}
\def\qq{\quad}
\let\ge\geqslant
\let\le\leqslant
\let\geq\geqslant
\let\leq\leqslant
\newcommand{\ca}{\begin{cases}}
\newcommand{\ac}{\end{cases}}
\newcommand{\ma}{\begin{pmatrix}}
\newcommand{\am}{\end{pmatrix}}
\renewcommand{\[}{\begin{equation}}
\renewcommand{\]}{\end{equation}}
\def\eq{\begin{equation}}
\def\qe{\end{equation}}
\def\[{\begin{equation}}
\def\bu{\bullet}

\title[Invariants for Laplacians on periodic graphs]
{Invariants for Laplacians on periodic graphs}

\date{\today}
\author[Evgeny Korotyaev]{Evgeny Korotyaev}
\address{Department of Mathematical Analysis, Saint-Petersburg State University, Universitetskaya nab. 7/9, St. Petersburg, 199034, Russia,
\ korotyaev@gmail.com, \
e.korotyaev@spbu.ru,}
\author[Natalia Saburova]{Natalia Saburova}
\address{Department of Mathematical Analysis, Algebra and Geometry, Northern (Arctic) Federal University, Severnaya Dvina emb. 17, Arkhangelsk, 163002, Russia,
 \ n.saburova@gmail.com, \ n.saburova@narfu.ru}


\begin{abstract}
We consider a Laplacian on periodic discrete graphs. Its spectrum  consists of a finite number of bands. In a class of periodic 1-forms, i.e., functions defined on edges of the periodic graph, we introduce a subclass of minimal forms with a minimal number $\cI$ of edges in their supports on the period. We obtain a specific decomposition of the Laplacian into a direct integral in terms of minimal forms, where fiber Laplacians (matrices) have the minimal number $2\cI$ of coefficients depending on the quasimomentum and show that the number $\cI$ is an invariant of the periodic graph. Using this decomposition, we estimate the position of each band, the Lebesgue measure of the Laplacian spectrum and the effective masses at the bottom of the spectrum in terms of the invariant $\cI$ and the minimal forms. In addition, we consider an inverse problem: we determine necessary and sufficient conditions for matrices depending on the quasimomentum on a finite graph to be fiber Laplacians. Moreover, similar results for Schr\"odinger operators with periodic potentials are obtained.
\end{abstract}

\maketitle

\section {\lb{Sec1}Introduction}
\setcounter{equation}{0}

Over the last 20 years the number of papers about Laplace and
Schr\"odinger operators on periodic discrete graphs significantly
increased due to their applications to problems of physics and
chemistry (see \cite{Ha02}, \cite{NG04}). The Floquet-Bloch theory for discrete operators on periodic graphs was discussed in \cite{HS99a}, \cite{KSS98}, \cite{SS92}. Some spectral properties of the Laplace and
Schr\"odinger operators on periodic discrete graphs were studied in \cite{FLP17}, \cite{HN09}, \cite{HS04}, \cite{KS14}, \cite{LP08}. Discrete magnetic Schr\"odinger operators on periodic graphs were considered in \cite{HS99a}, \cite{HS99b}, \cite{KS17}. It is known that the spectrum of the operators consists of an absolutely continuous part
(a~union of a finite number of non-degenerate bands) and a finite
number of flat bands, i.e., eigenvalues  of infinite multiplicity.
One of the main problems in this context is to obtain information about the band and gap structure; in particular to estimate the positions and length of bands and gaps. For example, in the one-dimensional case it was done
for Hill operators, \cite{K00}, \cite{K06}, where the proof is based on the
detailed analysis of the corresponding band functions as functions
of the quasimomentum. For the case of multi-dimensional graphs one
needs to investigate band functions of the quasimomentum in more
detail. This is a difficult problem that requires a detailed analysis of the dependence of the fiber Laplacian on the quasimomentum. In
particular, it is useful to obtain a representation of fiber
Laplacians with the minimal number of coefficients depending on the
quasimomentum.

We consider Schr\"odinger operators with periodic potentials on
periodic discrete graphs. We describe our main goals:

$\bu$ to decompose the Schr\"odinger operators into a direct integral, where fiber operators (matrices) have the minimal number $2\cI$ of coefficients depending on the quasimomentum; to show that the number $\cI$ is an invariant of the periodic graph;

$\bu$ to determine a localization of bands in terms of eigenvalues of
Schr\"odinger operators on some auxiliary finite graphs;

$\bu$ to obtain a \emph{sharp} estimate of the Lebesgue measure of the spectrum in terms of the invariant $\cI$, i.e., an estimate becoming an identity for specific graphs;

$\bu$ to estimate the effective masses at the bottom of the Laplacian spectrum in terms of geometric parameters of the graphs;

$\bu$ to solve an inverse problem: to determine necessary and sufficient conditions for matrices depending on the quasimomentum on a finite graph to be fiber Laplacians.

The proof of our results is essentially based on the notion of \emph{minimal forms} on graphs. A minimal form is a periodic function defined on edges of the periodic graph with a minimal support on the period.

\subsection{Periodic graphs.} Let $\cG=(\cV,\cE)$ be a connected infinite graph, possibly  having loops and multiple edges and embedded into the space $\R^d$. Here $\cV$ is the set of its vertices and $\cE$ is the set of its unoriented edges. Considering each edge in $\cE$ to have two orientations, we introduce the set $\cA$ of all oriented edges. An edge starting at a vertex $u$ and ending at a vertex $v$ from $\cV$ will be denoted as the ordered pair $(u,v)\in\cA$ and is said to be \emph{incident} to the vertices. Let $\ul\be=(v,u)$ be the inverse edge of $\be=(u,v)\in\cA$. Vertices $u,v\in\cV$ will be called \emph{adjacent} and denoted by $u\sim v$, if $(u,v)\in \cA$. We define the degree ${\vk}_v$ of
the vertex $v\in\cV$ as the number of all edges from
$\cA$, starting at $v$. A sequence of directed edges $(\be_1,\be_2,\ldots,\be_n)$ is called a \emph{path} if the terminus of the edge $\be_s$ coincides with the origin of the edge $\be_{s+1}$ for all $s=1,\ldots,n-1$. If the terminus of $\be_n$ coincides with the origin of $\be_1$, the path is called a \emph{cycle}.

Let $\G$ be a lattice of rank $d$ in $\R^d$ with a basis $\A=\{a_1,\ldots,a_d\}$, i.e.,
$$
\G=\Big\{a : a=\sum_{s=1}^dn_sa_s, \; n_s\in\Z,\; s\in\N_d\Big\}, \qqq \N_d=\{1,\ldots,d\},
$$
and let
\[\lb{fuce}
\Omega=\Big\{x\in\R^d : x=\sum_{s=1}^dt_sa_s, \; 0\leq t_s<1,\; s\in\N_d\Big\}
\]
be the \emph{fundamental cell} of the lattice $\G$. We define the equivalence relation on $\R^d$:
$$
x\equiv y \; (\hspace{-4mm}\mod \G) \qq\Leftrightarrow\qq x-y\in\G \qqq
\forall\, x,y\in\R^d.
$$

We consider \emph{locally finite $\G$-periodic graphs} $\cG$, i.e., graphs satisfying the
following conditions:
\begin{itemize}
  \item[1)] $\cG=\cG+a$ for any $a\in\G$;
  \item[2)] the quotient graph  $G_*=\cG/\G$ is finite.
\end{itemize}
The basis $a_1,\ldots,a_d$ of the lattice $\G$ is called the {\it periods} of $\cG$. We also call the quotient graph $G_*=\cG/\G$ the \emph{fundamental graph}
of the periodic graph $\cG$. The fundamental graph $G_*$ is a graph on the $d$-dimensional torus $\R^d/\G$. The graph $G_*=(\cV_*,\cE_*)$ has the vertex set $\cV_*=\cV/\G$, the set $\cE_*=\cE/\G$ of unoriented edges and the set $\cA_*=\cA/\G$ of oriented edges which are finite.

\subsection{Coordinate forms} For each $x\in\R^d$ we introduce the vector $x_\A\in\R^d$ by
\[\lb{cola}
x_\A=(t_1,\ldots,t_d), \qqq \textrm{where} \qq x=\textstyle\sum\limits_{s=1}^dt_sa_s.
\]
In other words, $x_\A$ is the coordinate vector of $x$ with respect to the basis $\A=\{a_1,\ldots,a_d\}$ of the lattice $\G$.

On the set $\cA$ of all oriented edges of the $\G$-periodic graph $\cG$ we define the surjection
\[\lb{sur}
\gf:\cA\rightarrow\cA_*=\cA/\G,
\]
which maps each $\be\in\cA$ to its equivalence class $\be_*=\gf(\be)$ which is an oriented edge of the fundamental graph $G_*$.

For any oriented edge $\be=(u,v)\in\cA$ of the periodic graph $\cG$ the \emph{edge coordinates} $\k(\be)$ are defined as the vector in $\R^d$ given by
\[
\lb{edco}
\k(\be)=v_\A-u_\A\in\R^d.
\]

On the fundamental graph $G_*$ we introduce the vector-valued \emph{coordinate form} $\k: \cA_*\ra\R^d$:
\[
\lb{dco}
\k(\be_*)=\k(\be) \qq \textrm{for some \; $\be\in\cA$ \; such that }  \; \be_*=\gf(\be), \qqq \be_*\in\cA_*,
\]
where $\gf$ is defined by \er{sur}.

The coordinate form $\k$ is uniquely determined by \er{dco}, since
$$
\k(\be+a)=\k(\be),\qqq \forall\, (\be,a)\in\cA \ts\G.
$$

\subsection{Minimal forms} We introduce the notion of minimal forms on graphs which will be used in the formulation of our results. Let $G=(V,E)$ be a finite connected graph with the vertex set $V$, the set $E$ of unoriented edges and the set $A$ of oriented edges. Denote by $\cC$ the cycle space of the graph $G$.

A vector-valued function $\gb: A\ra\R^d$ satisfying the condition $\gb(\ul\be\,)=-\gb(\be)$ for all $\be\in A$ is called a \emph{1-form} on the graph $G$. For any 1-form $\gb$ we define the vector-valued \emph{flux function} $\Phi_\gb:\cC\to\R^d$ as follows:
\[\lb{mafla}
\Phi_\gb(\mathbf{c})=\sum_{\be\in\mathbf{c}}\gb(\be), \qqq \mathbf{c}\in\cC.
\]
We call $\Phi_\gb(\mathbf{c})\in\R^d$ the \emph{flux} of the 1-form $\gb$ through the cycle $\mathbf{c}$.

We fix a 1-form $\mathbf{x}: A\ra\R^d$ and define the set $\mF(\mathbf{x})$ of all 1-forms $\gb: A\ra\R^d$ satisfying the condition $\Phi_\gb=\Phi_{\mathbf{x}}$:
\[\lb{vv1f}
\mF(\mathbf{x})=\{\gb: A\ra\R^d : \gb(\ul\be\,)=-\gb(\be)\qq \textrm{for all} \qq \be\in A, \qq \textrm{and} \qq \Phi_\gb=\Phi_{\mathbf{x}}\},
\]
where $\Phi_{\mathbf{x}}$ is the flux function of the form $\mathbf{x}$ defined by \er{mafla}.

In the set $\mF(\mathbf{x})$ we specify a subset of \emph{minimal} 1-forms and a subset of \emph{maximal} 1-forms. Let $\#M$ denote the number of elements in a set $M$.

\medskip

$\bu$ A 1-form $\gM\in\mF(\mathbf{x})$ is called \emph{maximal} if
\[\lb{dmMf}
\#\supp \gM\geq\#\supp \gb \qqq \textrm{for any 1-form} \qqq \gb\in\mF(\mathbf{x}).
\]

$\bu$ A 1-form $\gm\in\mF(\mathbf{x})$ is called \emph{minimal} if
\[\lb{dmf}
\#\supp \gm\leq\#\supp \gb \qqq \textrm{for any 1-form} \qqq \gb\in\mF(\mathbf{x}).
\]
Among all 1-forms $\gb\in\mF(\mathbf{x})$ the maximal form $\gM$ has a maximal support and the minimal form $\gm$ has a minimal support. Since any 1-form $\gb\in\mF(\mathbf{x})$ has a finite support, minimal forms $\gm\in\mF(\mathbf{x})$ and maximal forms $\gM\in\mF(\mathbf{x})$ exist. For example, the coordinate form $\k$ defined by \er{edco}, \er{dco} is a maximal form in the set $\mF(\k)$ on the fundamental graph $G_*$ (see Proposition \ref{Pal0}.\emph{i}). An explicit expression for all minimal forms $\gm\in\mF(\mathbf{x})$ and for the number of edges in their supports are given in Theorem \ref{Pphi}.

\subsection{Schr\"odinger operators on periodic graphs.}

Let $\ell^2(\cV)$ be the Hilbert space of all square summable functions $f:\cV\to \C$ equipped
with the norm
$$
\|f\|^2_{\ell^2(\cV)}=\sum_{v\in\cV}|f(v)|^2<\infty.
$$

We define the \emph{discrete combinatorial Laplacian} $\D$ on
$f\in\ell^2(\cV)$ by
\[\lb{DLO}
\big(\D f\big)(v)=\sum_{(v,u)\in\cA}\big(f(v)-f(u)\big), \qqq
v\in\cV.
\]
The sum in \er{DLO} is taken over all oriented edges starting at the vertex $v$. It is well known (see, e.g., \cite{M91}) that $\D$
is self-adjoint and the point 0 belongs to its spectrum $\s(\D)$ containing in $[0,2\vk_+]$, i.e.,
\[
\lb{bf}
0\in\s(\D)\subset[0,2\vk_+],\qqq
\textrm{where}\qqq
\vk_+=\max_{v\in\cV}\vk_v<\iy.
\]

We consider the Schr\"odinger operator $H$ acting on the Hilbert space $\ell^2(\cV)$ and given by
\[
\lb{Sh}
H=\D+Q,
\]
where the potential $Q$ is real valued and satisfies
for all $(v,a)\in\cV\ts\G$:
\[
\lb{Pot} \big(Q f\big)(v)=Q(v)f(v), \qqq Q(v+a)=Q(v).
\]

\medskip

The paper is organized as follows. In Section \ref{Sec1'} we formulate our main results:

$\bu$ the direct integral decomposition for the Schr\"odinger operators $H$ in terms of minimal forms $\gm\in\mF(\k)$, where $\k$ is the coordinate form  on the fundamental graph $G_*$ defined by \er{edco}, \er{dco}, and the fiber Laplacian has the minimal number $2\cI=\#\supp \gm$ of coefficients depending on the quasimomentum; the number $\cI$ is an \emph{invariant} of the $\G$-periodic graph $\cG$, i.e., $\cI$ does not depend on the choice of the embedding of $\cG$ into $\R^d$, the minimal form $\gm$ and the basis $a_1,\ldots,a_d$ of the lattice $\G$ (Theorem \ref{TDImf});

$\bu$ estimates of the invariant $\cI$ (Proposition \ref{TNNI});

$\bu$ necessary and sufficient conditions for matrices depending on the quasimomentum on a finite graph to be fiber Laplacians (Corollary \ref{TCo3});

$\bu$ a localization of bands and a sharp estimate of the Lebesgue measure of the spectrum of the Schr\"odinger operators $H$ in terms of the invariant $\cI$ (Theorem \ref{T1});

$\bu$ estimates on the effective mass tensor at the bottom of the spectrum for the Laplacians $\D$ (Theorem \ref{TTsem}).

In Section \ref{Sec2} we give a full description of minimal forms on a finite connected graph (Theorem \ref{Pphi}), which are essentially used in the proof of our results.

Section \ref{Sec3} is devoted to the direct integral decomposition for Schr\"odinger operators on periodic graphs. In this section we consider some properties of the coordinate form $\k$ on the fundamental graph $G_*$ (Proposition \ref{Pal0}) and give a general representation of fiber Laplacians in terms of any 1-form $\gb\in\mF(\k)$ (Theorem \ref{TDIf}). In Section \ref{Sec3} we also prove Theorem \ref{TDImf}, Proposition \ref{TNNI} and Corollary \ref{TCo3}.

In Section \ref{Sec4}, using the representation of fiber Laplacians in terms of minimal forms, we prove Theorem \ref{T1} about spectral estimates for the Schr\"odinger operators and show that these estimates become identities for specific graphs (Proposition \ref{TG1}). In this section we also prove Theorem \ref{TTsem}.

\section {\lb{Sec1'}Main results}
\setcounter{equation}{0}

\subsection{Decomposition of Schr\"odinger operators}
Recall that in the standard case  the Laplacian $-\D$ in $L^2(\R^d)$
is unitarily equivalent to the constant fiber direct integral:
$$
U(-\D)U^{-1}=\int^\oplus_{\T^d}\D(\vt)d\vt, \qqq
\D(\vt)=(i\na+\vt)^2, \qqq \T^d=\R^d/(2\pi\Z)^d,
$$
where $U:L^2(\R^d)\to \int^\oplus_{\T^d}L^2(\Omega)d\vt$ is the Gelfand transformation. Each fiber operator $\D(\vt), \vt\in\T^d$,
is a magnetic Laplacian with the magnetic vector potential $\vt$,
on $L^2(\Omega)$. The parameter $\vt$ is called the
\emph{quasimomentum}.

In the case of Laplace and Schr\"odinger operators on periodic graphs we have the similar decomposition (see Theorem \ref{TFD1}). Here we present a specific direct integral for Schr\"odinger operators on periodic graphs, where fiber operators have the minimal number of coefficients depending on the quasimomentum.

We introduce the Hilbert space
\[\lb{Hisp}
\mH=L^2\Big(\T^{d},{d\vt\/(2\pi)^d}\,,\cH\Big)=\int_{\T^{d}}^{\os}\cH\,{d\vt
\/(2\pi)^d}\,, \qqq \cH=\ell^2(\cV_*),
\]
i.e., a constant fiber direct integral equipped with the norm
$$
\|g\|^2_{\mH}=\int_{\T^{d}}\|g(\vt,\cdot)\|_{\ell^2(\cV_*)}^2\frac{d\vt
}{(2\pi)^d}\,,
$$
where the function $g(\vt,\cdot)\in\cH$ for almost all
$\vt\in\T^{d}$.

\begin{theorem}
\label{TDImf}
Let $\gm\in\mF(\k)$ be a minimal form on the fundamental graph $G_*=(\cV_*,\cE_*)$, where $\mF(\k)$ is given by \er{vv1f} at $\mathbf{x}=\k$, and $\k$ is the coordinate form defined by \er{edco}, \er{dco}. Then

i) The image of the minimal form $\gm$ generates the group $\Z^d$. Moreover, for some basis of the lattice $\G$ there exist edges $\be_1,\ldots,\be_d\in\supp\gm$ such that the set
$\big\{\gm(\be_s)\big\}_{s\in\N_d}$ forms an orthonormal basis of $\Z^d$.

ii) The Schr\"odinger operator $H=\D+Q$ on $\ell^2(\cV)$ has the following decomposition into a constant fiber direct integral
\[
\lb{raz1m}
\mU_\gm H \mU_\gm^{-1}=\int^\oplus_{\T^{d}}H_\gm(\vt)\frac{d\vt
}{(2\pi)^d}\,,
\]
where the unitary operator $\mU_\gm:\ell^2(\cV)\to\mH$ is a product of the Gelfand transformation $U$ and the gauge transformation $W_\gm$ defined by \er{5001} and \er{Uvt}, respectively.  Here the fiber Schr\"odinger operator $H_\gm(\vt)$ and the fiber Laplacian $\D_\gm(\vt)$ are given by
\[
\label{Hvtm}
H_\gm(\vt)=\D_\gm(\vt)+Q,\qqq \forall\,\vt\in\T^{d},
\]
\[
\label{l2.13am}
\begin{aligned}
\big(\D_\gm(\vt)f\big)(v)=\vk_vf(v)
-\sum_{\be=(v,u)\in\cA_*}e^{i\lan\gm(\be),\,\vt\ran}f(u), \qqq f\in\ell^2(\cV_*),\qqq v\in\cV_*,
\end{aligned}
\]
where $\lan\cdot\,,\cdot\ran$ denotes the standard inner product in $\R^d$.

iii) Let $\gm'\in\mF(\k)$ be another minimal form. Then for each $\vt\in\T^{d}$ the operators $\D_\gm(\vt)$ and $\D_{\gm'}(\vt)$ are unitarily equivalent.

iv) The number of exponents $e^{i\lan\gm(\be),\,\cdot\,\ran}\neq1$, $\be\in\cA_*$, in the identities \er{l2.13am} for the fiber Laplacians $\D_\gm(\cdot)$ is equal to
$2\cI$, where
\[\lb{dIm}
\cI=\textstyle\frac12\,\#\supp \gm.
\]

v) The number $\cI$ is an \textbf{invariant} of the periodic graph $\cG$, i.e., it does not depend on the choice of

$\bu$ the minimal form $\gm\in\mF(\k)$;

$\bu$ the embedding of $\cG$ into $\R^d$;

$\bu$ the basis $a_1,\ldots,a_d$ of the lattice $\G$.
\end{theorem}

\no \textbf{Remarks.} 1) For simple periodic graphs (the $d$-dimensional lattice, the hexagonal lattice, the Kagome lattice, etc.) it is not difficult to find a minimal form $\gm\in\mF(\k)$ on their fundamental graphs. But for an arbitrary periodic graph this may be a rather complicated problem (see Theorem \ref{Pphi}).

2) We can consider the fiber Laplacian $\D_\gm(\vt)$, $\vt\in
\T^{d}$, as a discrete magnetic Laplacian with the
magnetic vector potential $\a(\be)=\lan\gm(\be),\vt\ran$, $\be\in\cA_*$, on
the fundamental graph $G_*$. Note that the operator $\D_\gm(0)$ is
just the Laplacian defined by \er{DLO} on $G_*$.

3) Decompositions of Schr\"odinger operators
on periodic discrete graphs into a constant fiber direct integral were given in \cite{KS14}, \cite{S13}. In \cite{KS14} the fiber operator was expressed in terms of a special 1-form $\t\in\mF(\k)$ defined by \er{in}, \er{inf} (see Theorem \ref{TFD1}). In \cite{S13} the fiber operator was given in terms of the coordinate form $\k$, which is a maximal form on the fundamental graph $G_*$.

4) In Theorem \ref{TDIf} we give a general decomposition of Schr\"odinger operators on periodic discrete graphs into a constant fiber direct integral, where fiber operators $\D_\gb(\vt)$ are expressed in terms of any 1-form $\gb\in\mF(\k)$. Moreover, we show that among all fiber Laplacians $\D_\gb(\vt)$, $\gb\in\mF(\k)$, the fiber operator $\D_\gm(\vt)$ has the minimal number of coefficients depending on the quasimomentum $\vt$.

5) The set $\mF(\k)$ is defined in terms of the coordinate form $\k$. But we can introduce this set using another its representative, for example, the 1-form $\t\in\mF(\k)$ defined by \er{in}, \er{inf}.

6) The minimal forms $\gm$, of course, depend on the choice of the embedding of the periodic graph $\cG$ into the space $\R^d$ and on the choice of the basis $a_1,\ldots,a_d$ of the lattice $\G$. But, due to Theorem \ref{TDImf}.\emph{v}) the number $2\cI$ of entries in $\#\supp \gm$ does not depend on these choices.

\medskip

We formulate some simple properties of the invariant $\cI$.

\begin{proposition}
\lb{TNNI} The invariant $\cI$ of the $\G$-periodic graph $\cG$ given by \er{dIm} satisfies
\[\lb{nni}
d\leq\cI\leq\b, \qqq \b=\#\cE_*-\#\cV_*+1,
\]
where $d$ is the rank of the lattice $\G$ and  $\b$ is the
Betti number of the fundamental graph $G_*=(\cV_*,\cE_*)$. Moreover,

i) there exists a periodic graph such that $d=\cI=\b$;

ii) there exists a periodic graph such that $d=\cI<\b$;

iii) there exists a periodic graph such that $d<\cI=\b$;

iv) there exists a periodic graph such that $d<\cI<\b$;

v) for any nonnegative integer $\mn$, there exists a periodic graph such that $\b-\cI=\mn$.
\end{proposition}

Now we consider the inverse problem: when a matrix-valued function $\mA(\vt)$ depending on the quasimomentum $\vt$ on a finite connected graph $G$ is a fiber Laplacian for some periodic graph $\cG$.

\begin{corollary}\label{TCo3}
Let $\gm: A\ra\R^d$ be a minimal form on a finite connected graph $G=(V,E)$ satisfying $\Phi_{\gm}(\cC)=\Z^d$, where $\cC$ is the cycle space of $G$, and $\Phi_{\gm}$ is the flux function of the form $\gm$ given by \er{mafla}. Then the operator $\mA(\vt):\ell^2(V)\to\ell^2(V)$, $\vt\in\T^d$, defined by
\[
\label{l2.13at}
\begin{aligned}
\big(\mA(\vt)f\big)(v)=\vk_vf(v)
-\sum_{\be=(v,u)\in A}e^{i\lan\gm (\be),\,\vt\ran}f(u), \qqq v\in V,
\end{aligned}
\]
is a fiber operator for the Laplacian $\D$ on some periodic graph $\cG$ with the fundamental graph $G_*=G$.
\end{corollary}

\no \textbf{Remark.} This corollary and Theorem \ref{TDImf}.\emph{ii}) give  necessary and sufficient conditions for an operator of the form \er{l2.13at} to be a fiber Laplacian on a given graph $G$.

\subsection{Spectrum of the Schr\"odinger operator.}
Theorem \ref{TDImf} and standard arguments (see Theorem XIII.85 in
\cite{RS78}) describe the spectrum of the Schr\"odinger operator
$H=\D+Q$. Each fiber operator $H_\gm(\vt)$, $\vt\in\T^d$,
has $\n$ eigenvalues $\l_{n}(\vt)$, $n\in\N_\n$, $\n=\#\cV_*$, which are labeled (counting multiplicities) by
\[
\label{eq.3} \l_{1}(\vt)\leq\l_{2}(\vt)\leq\ldots\leq\l_{\nu}(\vt),
\qqq \forall\,\vt\in\T^d.
\]
Since $H_\gm(\vt)$ is self-adjoint and analytic in $\vt\in\T^d$, each $\l_{n}(\cdot)$, $n\in\N_\n$, is a real and piecewise analytic function on the torus $\T^d$ and creates the \emph{spectral band} (or \emph{band} for short) $\s_n(H)$ given by
\[
\lb{ban.1}
\s_{n}(H)=[\l_{n}^-,\l_{n}^+]=\l_{n}(\T^{d}).
\]
Note that $\l_{1}^{-}=\l_{1}(0)$ (see \cite{SS92}).

\smallskip

Thus, the spectrum of the Schr\"odinger operator $H$ on the periodic graph $\cG$ is
given by
\[\lb{spec}
\s(H)=\bigcup_{\vt\in\T^d}\s\big(H_\gm(\vt)\big)=\bigcup_{n=1}^{\nu}\s_n(H).
\]
Note that if $\l_{n}(\cdot)= C_{n}=\mathrm{constant}$ on some subset of $\T^d$ of
positive Lebesgue measure, then  the operator $H$ on $\cG$ has the
eigenvalue $C_{n}$ of infinite multiplicity. We call $C_{n}$
a \emph{flat band}.

Thus, the spectrum of the Schr\"odinger operator
$H$ on the periodic graph $\cG$ has the form
$$
\s(H)=\s_{ac}(H)\cup \s_{fb}(H),
$$
where $\s_{ac}(H)$ is the absolutely continuous spectrum, which is a
union of non-degenerate bands, and $\s_{fb}(H)$ is the set of
all flat bands (eigenvalues of infinite multiplicity). An open
interval between two neighboring non-degenerate bands is
called a \emph{gap}.

\subsection{Estimates of the Lebesgue measure of the spectrum}
Now we estimate the position of the bands $\s_n(H)$, $n\in\N_\n$,
defined by \er{ban.1} in terms of eigenvalues of the Schr\"odinger operator on some subgraph of the fundamental graph $G_*$ and the maximum vertex degree of the remaining subgraph of $G_*$. We also obtain the estimate of the Lebesgue measure of the spectrum of the Schr\"odinger operator $H$ on a periodic graph $\cG$ in terms of the invariant $\cI$ given by \er{dIm}.

Let $\gb\in\mF(\k)$, where $\mF(\k)$ is given by \er{vv1f} at $\mathbf{x}=\k$, and $\k$ is the coordinate form defined by \er{edco}, \er{dco}. In order to formulate our result we equip each unoriented edge $\be\in\cE_*$ of the fundamental graph $G_*=(\cV_*,\cE_*)$ with some orientation and represent $G_*$ as a union of two graphs with the same vertex set $\cV_*$:
\[\lb{cEgm}
\begin{aligned}
G_*=G_\gb^0\cup G_\gb, \qqq G_\gb=(\cV_*,\cE_\gb),\qqq G_\gb^0=(\cV_*,\cE_*\sm\cE_\gb),\\ \textrm{where}\qqq \cE_\gb=\cE_*\cap\supp \gb.
\end{aligned}
\]
The graph $G_\gb^0$ is obtained from the fundamental graph $G_*$ by deleting all edges of $\supp \gb$ and preserving the vertex set $\cV_*$. Let $\D_\gb^0$ be the Laplacian defined by \er{DLO} on the graph $G_\gb^0$. Since the graph $G_\gb^0$ consists of $\n=\#\cV_*$ vertices, the Schr\"odinger operator $H_\gb^0=\D_\gb^0+Q$ on the graph $G_\gb^0$ has $\n$ eigenvalues $\m_{\gb,n}$, $n\in\N_\n$, which are labeled
(counting multiplicities) by
\[\lb{evH0}
\m_{\gb,1}\leq\ldots\leq\m_{\gb,\n}.
\]

\begin{theorem}
\lb{T1} i) Each band $\s_n(H)$, $n\in\N_\n$, of the Schr\"odinger operator $H=\D+Q$ on a periodic graph $\cG$
defined by \er{ban.1} satisfies
\[\lb{fesbp1}
\s_n(H)\subset\big[\m_{\gm,n},\,\m_{\gm,n}+2\vk^{\gm}_+\big],
\qq \textrm{where}\qq
\vk^{\gm}_+=\max\limits_{v\in\cV_*}\vk_v^{\gm},
\]
$\gm\in\mF(\k)$ is a minimal form, $\mu_{\mathfrak{m},n}$ are given by (\ref{evH0}) as $\mathfrak{b}=\mathfrak{m}$, and  $\vk_v^{\gm}$ is the degree of the vertex $v\in\cV_*$ on the graph $G_\gm$ defined in \er{cEgm}.

ii) The Lebesgue measure $|\s(H)|$ of the spectrum of the Schr\"odinger operator $H$ satisfies
\[
\lb{eq.7'}
|\s(H)|\le \sum_{n=1}^{\n}|\s_n(H)|\le 4\cI,
\]
where the invariant $\cI$ is given by \er{dIm}.

iii) Let $\mm,\mn\in\N$ and $\mm\leq\mn$. Then there exists a periodic graph $\cG$ such that $\cI=d=\mm$, $\b=\mn$ and the estimate \er{eq.7'} becomes an identity for any potential $Q$.
\end{theorem}

\no \textbf{Remarks.} 1) In \cite{FLP17} the authors also obtained a localization of the bands $\s_n(\D)$, $n\in\N_\n$, for the Laplacian $\D$, using the so-called \emph{index form} $\t\in\mF(\k)$ defined by \er{in}, \er{inf}. Let $\m_{\t,1}\leq\ldots\leq\m_{\t,\n}$ be the eigenvalues of the Laplacian $\D_\t^0$ on the graph $G_\t^0$ defined in \er{cEgm} as $\gb=\t$. Denote by $\cV_\t\ss\cV_*$ the set of fundamental graph vertices such that each edge $\be\in\cE_\t=\cE_*\cap\supp \t$ is incident to at least one vertex of the set $\cV_\t$. Let
\[
\lb{LDN} \m_{\t,1}^D\le \m_{\t,2}^D\le\ldots\le \m_{\t,\n-r}^D,\qqq
r=\# \cV_\t, \qqq r\geq1,
\]
be the eigenvalues of the Laplacian defined by \er{DLO} on $f\in \ell^2(\cV_*)$ with Dirichlet boundary condition $f|_{\cV_\t}=0$. Then each band $\s_n(\D)$, $n\in\N_\n$, satisfies \cite{FLP17}
\[\lb{lpl}
\begin{array}{ll}
\s_n(\D)\ss\big[\m_{\t,n},\m_{\t,n}^D\big],\qq & n=1,\ldots,\n-r, \\[6pt]
\s_n(\D)\ss\big[\m_{\t,n},2\vk_+\big],
\qqq & n=\n-r+1,\ldots,\n,
\end{array}
\]
where $\vk_+=\max\limits_{v\in\cV_*}\vk_v$, $\vk_v$ is the degree of the vertex $v\in\cV_*$ on the fundamental graph $G_*$. Finally the authors wrote "\emph{The bracketing intervals \er{lpl} depend on the fundamental domain. A good choice would be one where the set of connecting edges} (i.e., edges from $\supp \t$) \emph{is as small as possible}". For simple periodic graphs (the $d$-dimensional lattice, the triangular lattice, the hexagonal lattice, etc.) it is not difficult to find the index form $\t$ with the minimal number of edges in its support. However, in general, the number of edges in $\supp \t$ is difficult to control and it essentially depends on the choice of the embedding of the periodic graph $\cG$ into the space $\R^d$. Thus, in order to obtain a more precise localization than \er{lpl} one needs to take the minimal form $\gm\in\mF(\k)$ instead of the index form $\t$, since the number $\#\supp\t$ can be significantly greater than $\#\supp\gm$ which is an invariant of the periodic graph.

2) There exist other localizations of the bands of the Laplace and Schr\"odinger operators on periodic discrete graphs (see the results overview below). These localizations can be combined to obtain a more exact estimate of spectral band positions.

3) For most periodic graphs $|\s(H)|<4\cI$, i.e., the inequality \er{eq.7'} is strict. For example, the spectrum of the Laplacian $\D$ on the Kagome lattice (see Fig. \ref{ff.0.11}) has the form
$$
\s(\D)=\s_{ac}(\D)\cup\s_{fb}(\D),\qqq \s_{ac}(\D)=[0,6],\qqq \s_{fb}(\D)=\{6\}.
$$
The invariant $\cI=\textstyle\frac12\,\#\supp \gm$ is equal to 3 (one of the minimal forms $\gm$ is shown in Fig. \ref{ff.0.11}\emph{b}). Thus, for the Kagome lattice the estimate \er{eq.7'} has the form $|\s(\D)|=6<12=4\cdot\cI$.

4) It is known that $|\s(\D)|\leq2\vk_+$, where $\vk_+$ is defined in \er{bf}.  Thus, for the Laplacian $\D$ the estimate \er{eq.7'} is effective if and only if $2\cI<\vk_+$. A class of periodic graphs for which the inequality $2\cI<\vk_+$ holds true is considered in Proposition \ref{TG1}. Moreover, the difference $\vk_+-2\cI$ can be arbitrarily large natural number (see the third formula in \er{deco}).

5) The decomposition of the Schr\"odinger operator $H$ into the direct integral \er{raz1m} -- \er{l2.13am} in terms of minimal forms is essentially used in the proof of the estimate \er{eq.7'}.

6) From Theorem \ref{T1}.\emph{iii}) it follows that
the Lebesgue measure of the spectrum can be arbitrarily large (for specific graphs).

\medskip

We know only a few papers about estimates on bands of Laplace and Schr\"odinger operators on periodic graphs.
Lled\'o and Post \cite{LP08} estimated the positions of bands of  Laplacians both on metric and discrete graphs in terms of eigenvalues of the operator on finite graphs (the so-called eigenvalue bracketing). In \cite{FLP17} the authors extended the previous article and presented a band localization \er{lpl} for the discrete Laplacians. They also obtained a simple geometric condition that guarantees the existence of gaps for the operators. Korotyaev and Saburova \cite{KS15}
described a localization of bands and estimated the Lebesgue measure of the spectrum of the Schr\"odinger operators $H$ with periodic potentials on periodic discrete graphs in terms of eigenvalues of Dirichlet and Neumann operators on a fundamental domain of the periodic graph. In \cite{KS14} they also estimated the Lebesgue measure of the spectrum of $H$ in terms of geometric parameters of the
graph:
\[
\lb{esLm0}
|\s(H)|\le \sum_{n=1}^{\n}|\s_n(H)|\le 4b, \qqq \textrm{where} \qqq b=\textstyle\frac12\,\#\supp\t,
\]
$\t\in\mF(\k)$ is the index form defined by \er{in}, \er{inf}. But in contrast to the invariant $\cI$ in \er{eq.7'}, this number $b$ depends essentially on the choice of the embedding of the periodic graph $\cG$ into the space $\R^d$, i.e., $b$ is not an invariant for $\cG$. Moreover, the number $b$ can be less, equal or significantly greater than the Betti number $\b$ of the fundamental graph $G_*$. For the magnetic Schr\"odinger operators $H_\a$ with a periodic magnetic potential $\a$ we obtained the following estimate (see \cite{KS17})
\[
\lb{esLm1}
|\s(H_\a)|\le \sum_{n=1}^{\n}|\s_n(H_\a)|\le 4\b.
\]
Due to Proposition \ref{TNNI}.\emph{v}) the difference between the Betti number $\b$ and the invariant $\cI$ given by \er{dIm} can be any natural number (for specific graphs). Then comparing estimates \er{eq.7'} and  \er{esLm1} as $\a=0$, we see that the estimate \er{eq.7'} is better than \er{esLm1}.

\subsection{Effective masses at the bottom of the spectrum for discrete Laplacians}
We introduce the notion of effective masses. Firstly, for simplicity, we
consider an isotropic case. For a free electron we have the energy
dispersion:
\[\lb{drfe}
\l_0(\vt)=\frac{\vt^2}{2m_0}\,,
\]
where $\vt\in\R^d$ is the momentum, $m_0$ is the mass of the free
electron. Now we consider an electron in a crystal lattice.
Typically the first band function $\l(\vt)$ of the Hamiltonian, where
$\vt\in \T^d=\R^d/(2\pi\Z)^d$ is the quasimomentum, has its minimum
at the point $\vt_0=0$. Then the function $\l(\vt)$ can be locally
approximated as
\[
\lb{tem'}
\begin{aligned}
\l(\vt)=\frac{\vt^2}{2m} +O(|\vt|^3).
\end{aligned}
\]
Here the linear terms vanish, since $\l(\vt)$ has an extremum at the
point $\vt_0=0$. Comparing dispersion relationships \er{drfe} and
\er{tem'}, we see that for energies close to the point $\vt_0=0$ the
electron can be treated as if it were free, but with a different mass
$m$. This mass $m$ is called \emph{the effective mass}. The effective mass approximation is a standard approach in physics when the complicated Hamiltonian, for small energies, is replaced by the model
Hamiltonian $-{\D_\vt\/2m}$, where $\D_\vt$ is the Laplacian in the quasimomentum space.

Secondly, we consider a general (anisotropic) case. Let $\l(\vt)$, $\vt\in\T^d$, be the first band function of the Laplacian $\D$. We recall that $\l(\cdot)$ is a piecewise analytic function. It is known that
\[\lb{bes}
\min\limits_{\vt\in\T^d}\l(\vt)=\l(0)=0,
\]
and $\l(0)=0$ is a simple eigenvalue of the operator $\D(0)$. Then the band function $\l(\cdot)$ has a Taylor series expansion as $|\vt|\to0$:
\[
\label{lam}
\l(\vt)={1\/2}\sum\limits_{i,j=1}^d\,M_{ij}\,\vt_i \vt_j+O\big(|\vt|^3\big),\qqq
M_{ij}={\pa^2\l(0)\/\pa\vt_i\pa\vt_j}\,,\qqq \vt=(\vt_j)_{j=1}^d.
\]
The matrix $m=M^{-1}$ represents \emph{the effective mass tensor} at the bottom of the spectrum, where $M=(M_{ij})_{i,j=1}^d$. In this general case the effective mass is a tensor and its components depend on the coordinate system in the quasimomentum space.

In \cite{K08} estimates on effective masses for zigzag
nanotubes in magnetic fields were obtained. In \cite{KS16} the authors estimated effective masses associated with the ends of each band for the Laplacians on periodic graphs in terms of geometric parameters of the graphs.

Using the direct integral decomposition \er{raz1m} -- \er{l2.13am} for Laplacians in terms of minimal forms we improve the estimates of the effective masses at the bottom of the spectrum obtained in \cite{KS16}. Estimates for the effective mass tensor are very important to study the discrete spectrum in a gap of the Schr\"odinger operator with a periodic potential perturbed by a decaying potential on periodic graphs, see \cite{KoS18}.

\begin{theorem}\label{TTsem}
Let $M=(M_{ij})_{i,j=1}^d$ be defined in \er{lam}. Then the effective mass tensor $m=M^{-1}$ at the bottom of the spectrum of the Laplacian $\D$ satisfies (in the sense of quadratic forms)
\[
\lb{Mtem}
{\nu\/C_\gm}\le m\le \frac{\n^2d}2\,, \qqq \textrm{where} \qqq C_\gm=\sum\limits_{\be\in\supp\gm}\|\gm(\be)\|^2,
\]
$\gm\in\mF(\k)$ is a minimal form on the fundamental graph $G_*=(\cV_*,\cE_*)$, $\n=\#\cV_*$. Moreover, if $d=\cI$, where $\cI=\textstyle\frac12\,\#\supp \gm$, then ${\n\/2d}\leq m\le\frac{\n^2d}2$.
\end{theorem}

\no {\bf Remarks.} 1) The upper estimate on the effective mass tensor $m$
in \er{Mtem} for some particular class of periodic graphs was proved in \cite{KS16}. Now we prove that this estimate holds true for any periodic graph. The lower estimate on the effective mass tensor at the bottom of the Laplacian spectrum obtained in \cite{KS16} has the form
\[
\lb{tseN'}
{\nu\/C_\t}\le m, \qqq \textrm{where} \qqq C_\t=\textstyle\sum\limits_{\be\in\supp\t}\|\t(\be)\|^2,
\]
$\t\in\mF(\k)$ is the index form defined by \er{in}, \er{inf}. In contrast to the invariant $2\cI=\#\supp \gm$, the number of entries of $\supp\t$ depends essentially on the choice of the embedding of the periodic graph $\cG$ into the space $\R^d$, i.e., this number is not an invariant for $\cG$.
Moreover, this number $\#\supp\t$ can be significantly greater than $2\cI$.

2) The identity $d=\cI$ holds true, for example, for the $d$-dimensional
lattice and the hexagonal lattice. For other examples of such graphs see Proposition \ref{TG1}.

\section{\lb{Sec2} Minimal forms and their properties}
\setcounter{equation}{0}

\subsection{Betti numbers and spanning trees}
We recall the definitions of the Betti number and spanning trees which will be used in the proof of our results. Let $G=(V,E)$ be a finite connected graph. Recall that $\#M$ denotes the number of elements in a set $M$.

\medskip

$\bullet$  A \emph{spanning tree} $T=(V,E_{T})$ of the graph $G$ is a connected subgraph of $G$ which has no cycles and contains all vertices of $G$.

$\bullet$ The \emph{Betti number} $\b$ of the graph $G$ is defined as
\[\lb{benu}
\b=\# E-\# V+1.
\]
Note that the Betti number $\b$ can also be defined in one of the following equivalent ways:
\begin{itemize}
  \item[\emph{i})] as the number of edges that have to be removed from $E$ to turn $G$ into a spanning tree of $G$;
  \item[\emph{ii})] as the dimension of the cycle space $\cC$ of the graph $G$, i.e.,
\[
\lb{debn}
\b=\dim\cC.
\]
\end{itemize}

\medskip

We introduce the set $\cS_T$ of all edges from $E$ that do not belong to the spanning tree $T$, i.e.,
\[
\lb{ulS}
\cS_T=E\sm E_T,
\]
and recall some properties of spanning trees (see, e.g., Lemma 5.1 and Theorem 5.2 in \cite{B74}, and Proposition 1.3 in \cite{CDS95}).

\medskip

\no\textbf{Properties of spanning trees.}\lb{PSTs}

\emph{1) The set $\cS_T$ consists of exactly $\b$ edges, where $\b$ is the Betti number defined by \er{benu}.}

\emph{2) A spanning tree $T$ is a maximally acyclic subgraph of $G$ and for each $\be\in\cS_T$ there exists a unique cycle $\mathbf{c}_\be$ on $G$
whose edges are all in $T$ except $\be$.}

\emph{3) The set of all such cycles $\mathbf{c}_\be$, $\be\in\cS_T$, forms a basis $\cB_T$ of the cycle space $\cC$ of the graph~$G$:}
\[\lb{bst}
\cB_T=\{\mathbf{c}_\be: \be\in\cS_T\}.
\]

\emph{4) Let $0=\x_1<\x_2\leq\x_3\leq\ldots\leq\x_\n$ be the eigenvalues of the Laplacian $\D$ on the graph $G=(V,E)$, where $\n=\# V$. Then the number of spanning trees of $G$ is equal to $\frac1\n\,\x_2\,\x_3\ldots\x_\n$.}

\medskip

\no \textbf{Example.} For the graph $G$ shown in Fig.\ref{ffS'}\emph{a} we can choose the spanning trees $T$ and $\wt T$ (Fig.\ref{ffS'} \emph{b,c}).
The set $\cS_T$ consists of three edges $\be_1,\be_2,\be_3$ (they are shown in Fig.\ref{ffS'} \emph{b,c} by the dotted lines) and depends on the choice of the spanning tree. The Betti number $\b$ defined by \er{benu} is equal to 3 and does not depend on the set $\cS_T$.

\setlength{\unitlength}{1.0mm}
\begin{figure}[h]
\centering
\unitlength 1.0mm 
\linethickness{0.4pt}
\begin{picture}(100,30)

\put(10,10){\circle{1}}
\put(10,30){\circle{1}}
\put(10,20){\circle{1}}
\put(0,20){\circle{1}}
\put(20,20){\circle{1}}
\put(10,10){\line(1,1){10.00}}
\put(10,10){\line(-1,1){10.00}}
\put(10,30){\line(1,-1){10.00}}
\put(10,30){\line(-1,-1){10.00}}

\put(10,30){\line(0,-1){10.00}}
\put(0,20){\line(1,0){20.00}}

\put(50,10){\circle{1}}
\put(50,30){\circle{1}}
\put(50,20){\circle{1}}
\put(40,20){\circle{1}}
\put(60,20){\circle{1}}
\put(50,10){\line(1,1){10.00}}
\put(50,10.2){\line(1,1){10.00}}
\put(50,9.8){\line(1,1){10.00}}
\put(50,10.1){\line(1,1){10.00}}
\put(50,9.9){\line(1,1){10.00}}

\put(50,30){\line(0,-1){10.00}}
\put(50.1,30){\line(0,-1){10.00}}
\put(49.9,30){\line(0,-1){10.00}}
\put(40,20.1){\line(1,0){20.00}}
\put(40,19.9){\line(1,0){20.00}}
\put(40,20){\line(1,0){20.00}}

\qbezier[20](50,10)(45,15)(40,20)
\qbezier[20](50,30)(55,25)(60,20)
\qbezier[20](50,30)(45,25)(40,20)

\put(42,13){$\be_1$}
\put(41,26){$\be_2$}
\put(55,26){$\be_3$}
\put(90,10){\circle{1}}
\put(90,30){\circle{1}}
\put(90,20){\circle{1}}
\put(80,20){\circle{1}}
\put(100,20){\circle{1}}
\put(90,10){\line(1,1){10.00}}
\put(90,10.2){\line(1,1){10.00}}
\put(90,9.8){\line(1,1){10.00}}
\put(90,10.1){\line(1,1){10.00}}
\put(90,9.9){\line(1,1){10.00}}

\put(90,30){\line(0,-1){10.00}}
\put(90.1,30){\line(0,-1){10.00}}
\put(89.9,30){\line(0,-1){10.00}}
\put(90,20.1){\line(1,0){10.00}}
\put(90,19.9){\line(1,0){10.00}}
\put(90,20){\line(1,0){10.00}}

\qbezier[15](80,20)(85,20)(90,20)
\qbezier[20](90,10)(85,15)(80,20)
\qbezier[20](90,30)(95,25)(100,20)

\put(90,30.2){\line(-1,-1){10.00}}
\put(90,29.8){\line(-1,-1){10.00}}
\put(90,30.1){\line(-1,-1){10.00}}
\put(90,29.9){\line(-1,-1){10.00}}
\put(90,30){\line(-1,-1){10.00}}
\put(82,13){$\be_1$}
\put(84,21){$\be_2$}
\put(95,26){$\be_3$}

\put(-3,10){(\emph{a})}
\put(33,10){(\emph{b})}
\put(73,10){(\emph{c})}

\put(0,26){$G$}
\put(33,26){$T$}
\put(80,26){$\wt T$}
\end{picture}

\vspace{-10mm}
\caption{\footnotesize  \emph{a}) A finite connected graph $G$;\quad \emph{b}),\,\emph{c}) the spanning trees $T$ and $\wt T$, $\cS_T=\{\be_1,\be_2,\be_3\}$, $\b=3$.} \label{ffS'}
\end{figure}
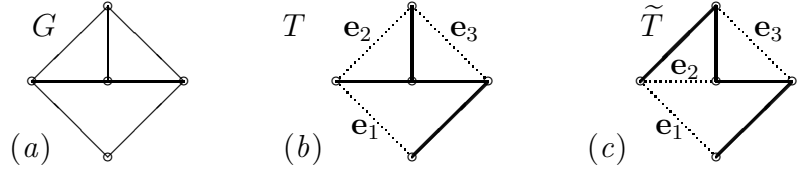

\subsection{Minimal forms}
Let $\mathbf{x}$ be a 1-form on a finite connected graph $G=(V,E)$. We describe all minimal forms $\gm\in\mF(\mathbf{x})$ on $G$, where $\mF(\mathbf{x})$ is defined by \er{vv1f}, and give an explicit expression for the number of edges in their supports.

Let $T=(V,E_{T})$ be a spanning tree of the graph $G$. Then, due to the properties 2) and 3) of spanning trees (see page \pageref{PSTs}), for each $\be\in\cS_T=E\sm E_T$, there exists a unique cycle $\mathbf{c}_\be$ consisting of only $\be$ and edges of $T$, and the set $\cB_T$ of all these cycles forms a basis of the cycle space $\cC$ of the graph $G$.

Among all spanning trees of the graph $G$ we specify a subset of \emph{minimal} spanning trees. Let $\b_T=\b_T(\mathbf{x})$ be the number of basic cycles from $\cB_T$ with non-zero fluxes of the form $\mathbf{x}$:
\[\lb{ncIT}
\b_T=\#\{\mathbf{c}\in\cB_T : \Phi_{\mathbf{x}}(\mathbf{c})\neq0\}.
\]
A spanning tree $T$ of the graph $G$ is called \emph{$\mathbf{x}$-minimal} if
\[\lb{dmf}
\b_T\leq\b_{\wt T} \qqq \textrm{for any spanning tree $\wt T$ of $G$.}
\]
Among all spanning trees of $G$ the $\mathbf{x}$-minimal one gives a minimal number of basic cycles from $\cB_T$ with non-zero fluxes of the form $\mathbf{x}$. Since the set $\cB_T$ of basic cycles is finite, minimal spanning trees exist. All $\mathbf{x}$-minimal spanning trees $T$ have the same number $\b_T$. We denote this number by $\b(\mathbf{x})$:
\[\lb{nbo}
\b(\mathbf{x})=\b_T \qqq \textrm{for any $\mathbf{x}$-minimal spanning tree $T$ of $G$.}
\]

\no \textbf{Example.} We consider the graph $G$ shown in Fig.\ref{FkfKl}\emph{a}. The values of the 1-form $\mathbf{x}:A\to\R^2$ on $G$ are shown in the figure near the corresponding edges. For the spanning tree $T$ shown in Fig.\ref{FkfKl}\emph{b} the set $\cS_T$ consists of 4 edges $\be_2,\be_4,\be_5,\be_6$ (they are shown by the dotted lines) and the set $\cB_T$ defined by \er{bst} consists of 4 basic cycles. Only one basic cycle $(\be_1,\be_2,\be_3)$ has zero flux of the form $\mathbf{x}$. Thus, the number $\b_T$ defined by \er{ncIT} is equal to $4-1=3$. For the spanning tree $\wt T$ (Fig.\ref{FkfKl} \emph{c}) all basic cycles from $\cB_{\wt T}$ have non-zero fluxes of the form $\mathbf{x}$. Thus, $\b_{\wt T}=4$. Looking over all spanning trees of $G$ (the number of spanning trees is finite) one can check that $T$ is a $\mathbf{x}$-minimal spanning tree of the graph $G$.

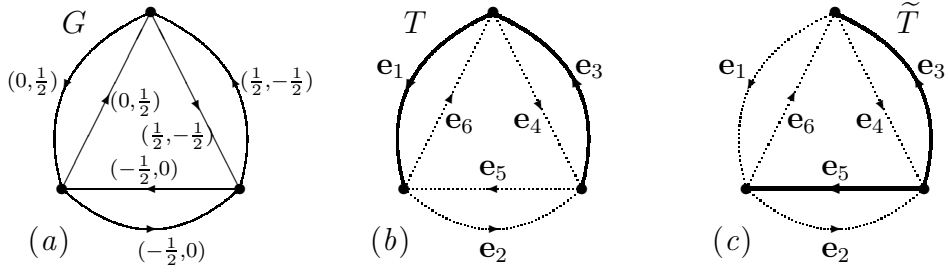
\begin{figure}[h]
\centering
\unitlength 0.9mm 
\linethickness{0.4pt}
\ifx\plotpoint\undefined\newsavebox{\plotpoint}\fi 
\hspace{-2.5cm}
\begin{picture}(146,48)(0,0)

\put(20,6){(\emph{a})}
\put(25,38){$G$}
\put(31.0,27.0){\vector(1,2){1.0}}
\put(44.2,28.5){\vector(1,-2){1.0}}
\put(38,15.0){\vector(-1,0){1.0}}

\put(25,15){\line(1,0){26.0}}
\put(25,15){\line(1,2){13.0}}

\put(50.8,29.5){\vector(-1,2){1.0}}
\put(26.2,31.5){\vector(-1,-2){1.0}}
\put(38,9.0){\vector(1,0){1.0}}

\put(51,15){\line(-1,2){13.0}}

\put(25,15){\circle*{1.5}}
\put(51,15){\circle*{1.5}}
\put(38,41){\circle*{1.5}}
\bezier{200}(25,15)(38,3)(51,15)
\bezier{200}(25,15)(20,33)(38,41)
\bezier{200}(51,15)(56,33)(38,41)

\put(32.0,17.0){$\scriptstyle(-{1\/2},0)$}
\put(32.0,27.0){$\scriptstyle(0,{1\/2})$}
\put(36.5,22.0){$\scriptstyle({1\/2},-{1\/2})$}

\put(36.0,5.0){$\scriptstyle(-{1\/2},0)$}
\put(17.0,30.0){$\scriptstyle(0,{1\/2})$}
\put(51.0,30.0){$\scriptstyle({1\/2},-{1\/2})$}

\put(70,6){(\emph{b})}
\put(75,38){$T$}
\bezier{40}(75,15)(88,15)(101,15)
\put(88,15.0){\vector(-1,0){1.0}}

\bezier{40}(75,15)(81.5,28)(88,41)
\put(81.2,27.2){\vector(1,2){1.0}}

\put(94.2,28.5){\vector(1,-2){1.0}}

\put(100.8,29.5){\vector(-1,2){1.0}}
\put(100.5,29.5){\vector(-1,2){1.0}}

\put(76.2,31.5){\vector(-1,-2){1.0}}
\put(76.5,31.5){\vector(-1,-2){1.0}}

\put(88,9.0){\vector(1,0){1.0}}

\bezier{40}(101,15)(94.5,28)(88,41)

\bezier{50}(75,15)(88,3)(101,15)

\bezier{200}(75,15)(70,33)(88,41)
\bezier{200}(75.2,15)(70.2,33)(88.2,40.8)

\bezier{200}(101,15)(106,33)(88,41)
\bezier{200}(100.8,15)(105.8,33)(87.8,40.7)

\put(75,15){\circle*{1.5}}
\put(101,15){\circle*{1.5}}
\put(88,41){\circle*{1.5}}

\put(86.0,17.0){$\be_5$}
\put(81.0,24.0){$\be_6$}
\put(91.0,24.0){$\be_4$}
\put(86.0,5.0){$\be_2$}
\put(71.0,32.0){$\be_1$}
\put(100.0,32.0){$\be_3$}

\put(147,38){$\wt T$}
\put(120,6){(\emph{c})}
\put(125,14.8){\line(1,0){26.0}}
\put(125,15){\line(1,0){26.0}}
\put(125,15.2){\line(1,0){26.0}}

\put(138,15.0){\vector(-1,0){1.0}}
\put(138,14.8){\vector(-1,0){1.0}}
\put(138,15.2){\vector(-1,0){1.0}}

\bezier{40}(125,15)(131.5,28)(138,41)

\put(131.2,27.5){\vector(1,2){1.0}}

\put(144.2,28.5){\vector(1,-2){1.0}}

\put(150.8,29.5){\vector(-1,2){1.0}}
\put(150.6,29.5){\vector(-1,2){1.0}}

\put(126.2,31.5){\vector(-1,-2){1.0}}

\put(138,9.0){\vector(1,0){1.0}}

\bezier{40}(151,15)(144.5,28)(138,41)

\bezier{50}(125,15)(138,3)(151,15)
\bezier{50}(125,15)(120,33)(138,41)

\bezier{200}(151,15)(156,33)(138,41)
\bezier{200}(150.8,15)(155.8,33)(137.8,40.7)

\put(125,15){\circle*{1.5}}
\put(151,15){\circle*{1.5}}
\put(138,41){\circle*{1.5}}

\put(136.0,17.0){$\be_5$}
\put(131.0,24.0){$\be_6$}
\put(141.0,24.0){$\be_4$}
\put(136.0,5.0){$\be_2$}
\put(121.0,32.0){$\be_1$}
\put(150.0,32.0){$\be_3$}

\end{picture}

\caption{\footnotesize  \emph{a}) A graph $G$ with the 1-form $\mathbf{x}$;\quad \emph{b}) a $\mathbf{x}$-minimal spanning tree $T$;\quad \emph{c}) a non-minimal spanning tree $\wt T$.}
\label{FkfKl}
\end{figure}

\medskip

Let $T=(V,E_T)$ be a spanning tree of the graph $G$. We equip each edge of $\cS_T=E\sm E_T$ with some orientation and associate with $T$ a 1-form $\gm(\,\cdot\,,T): A\to\R^d$ equal to zero on edges of $T$ and coinciding with the flux $\Phi_{\mathbf{x}}(\mathbf{c}_\be)$ of the form $\mathbf{x}$ through the cycle $\mathbf{c}_\be$ on each edge $\be\in\cS_T$, i.e.,
\[\lb{cat}
\gm(\be,T)=\left\{
\begin{array}{cl}
 \Phi_{\mathbf{x}}(\mathbf{c}_\be), &  \textrm{ if } \, \be\in\cS_T\\[6pt]
 -\Phi_{\mathbf{x}}(\mathbf{c}_{\ul\be}), &  \textrm{ if } \, \ul\be\in\cS_T\\[6pt]
  0, \qq & \textrm{ otherwise} \\
\end{array}\right..
\]

\begin{lemma}
\lb{L1f}
For each spanning tree $T=(V,E_{T})$ of the graph $G=(V,E)$ the function $\gm(\,\cdot\,,T): A\ra\R^d$ defined by \er{cat} is a 1-form in the set $\mF(\mathbf{x})$ given by \er{vv1f}.
\end{lemma}
\no \textbf{Proof.} From \er{cat} it follows that $\gm(\be,T)=-\gm(\ul\be,T)$ for each $\be\in A$. Thus, the function $\gm(\,\cdot\,,T)$ is a 1-form on $G$. For each $\be\in\cS_T=E\sm E_T$ there exists a unique cycle $\mathbf{c}_\be$ consisting of only $\be$ and edges of the spanning tree $T$. Then, due to \er{mafla}, \er{cat},
$$
\Phi_{\gm(\,\cdot\,,T)}(\mathbf{c}_\be)=\sum_{\wt\be\in\mathbf{c}_\be}\gm(\wt\be,T)=
\gm(\be,T)=\Phi_{\mathbf{x}}(\mathbf{c}_\be)
$$
for each basic cycle $\mathbf{c}_\be$, $\be\in\cS_T$. Consequently, for each cycle $\mathbf{c}\in\cC$ we have $\Phi_{\gm(\,\cdot\,,T)}(\mathbf{c})=\Phi_{\mathbf{x}}(\mathbf{c})$. Thus, $\gm(\,\cdot\,,T)\in\mF(\mathbf{x})$. \qq $\BBox$

\begin{theorem}
\lb{Pphi} Let $\mathbf{x}$ be a 1-form on a finite connected graph $G=(V,E)$, and $\mF(\mathbf{x})$ be given by \er{vv1f}. Then

i) Each minimal form $\wt\gm\in\mF(\mathbf{x})$ satisfies $\wt\gm=\gm(\,\cdot\,,T)$, where the 1-form $\gm(\,\cdot\,,T): A\ra\R^d$ is defined by \er{cat} for some $\mathbf{x}$-minimal spanning tree $T$ of the graph $G$.

ii) Conversely, for each $\mathbf{x}$-minimal spanning tree $T$ of the graph $G$ the 1-form $\gm(\,\cdot\,,T)$ is a minimal form in the set $\mF(\mathbf{x})$.

iii) The number $\b(\mathbf{x})$ defined by \er{ncIT}, \er{nbo}
satisfies
\[\lb{gfe1}
\b(\mathbf{x})=\textstyle\frac12\,\#\supp \gm,
\]
where $\gm\in\mF(\mathbf{x})$ is a minimal form.

iv) There exist a finite connected graph $\wt G$ and a 1-form $\wt{\mathbf{x}}$ on it such that $\gm(\,\cdot\,,T)=\gm(\,\cdot\,,\wt T\,)$ for some distinct $\wt{\mathbf{x}}$-minimal spanning trees $T$ and $\wt T$ of the graph $\wt G$.
\end{theorem}

\no \textbf{Proof.} \emph{i}) Let $\wt\gm\in\mF(\mathbf{x})$ be a minimal form on the graph $G=(V,E)$.
We consider the finite graph
\[\lb{Ggb}
G_0=(V,E_0), \qqq \textrm{where} \qqq E_0=E\sm\supp\wt\gm,
\]
and show that $G_0$ is connected. Let $c_0$ be the number of connected components of the graph $G_0$ and $\b_0$ be the Betti number of $G_0$:
\[\lb{bno}
\b_0=\# E_0-\#V+c_0.
\]
Since the graph $G$ is connected, then there exist $c_0-1$ edges $\be_1,\ldots,\be_{c_0-1}\in (E\cap\supp\wt\gm)$ such that the graph
\[\lb{Ggb1}
G_1=(V,E_1), \qqq \textrm{where} \qqq E_1=E_0\cup\{\be_1,\ldots,\be_{c_0-1}\},
\]
is connected. Let $T=(V,E_T)$ be a spanning tree of the graph $G_1$. This tree $T$ is also a spanning tree of $G$.  We consider a 1-form $\gm(\,\cdot\,,T)\in\mF(\mathbf{x})$ defined by \er{cat} on the graph $G$.
Since all $\b_0$ basic cycles of the graph $G_0$ have zero fluxes of the 1-form $\mathbf{x}$, using \er{cat}, \er{bno} and the identity $\#\cS_T=\b$, we obtain
\[\lb{smT}
\begin{aligned}
&\textstyle\frac12\,\#\supp\gm(\,\cdot\,,T)\leq\b-\b_0=
\#E-\# V+1-(\# E_0-\# V+c_0)=\\&\textstyle\# E-\# V+1-(\# E-\frac12\,\#\supp\wt\gm-\# V+c_0)=\frac12\,\#\supp\wt\gm+1-c_0.
\end{aligned}
\]
Thus, we have
\[\lb{smT1}
\textstyle\frac12\,\#\supp\gm(\,\cdot\,,T)+c_0-1\leq\frac12\,\#\supp\wt\gm.
\]
Since $\wt\gm$ is a minimal form in the set $\mF(\mathbf{x})$ and $\gm(\,\cdot\,,T)\in\mF(\mathbf{x})$, we obtain that $\#\supp\gm(\,\cdot\,,T)=\#\supp\wt\gm$ and $c_0=1$. The first identity and the definition \er{cat} of the form $\gm(\,\cdot\,,T)$ yield
\[\lb{sugm}
\textstyle\frac12\#\supp\wt\gm=\frac12\#\supp\gm(\,\cdot\,,T)=\b_T,
\]
where $\b_T$ is defined by \er{ncIT}. From the identity $c_0=1$ it follows that the graph $G_0$ defined by \er{Ggb} is connected, $G_1=G_0$, and $T=(V,E_T)$ is a spanning tree of $G_0$.

We show that $\wt\gm=\gm(\,\cdot\,,T)$. Let $\be\in\cS_T$, where $\cS_T=E\sm E_{T}$. Then on $G$ there exists a unique cycle $\mathbf{c}_\be$ whose edges are all in $T$ except $\be$. Due to the definition \er{Ggb} of the graph $G_0$,  $\supp\wt\gm\cap E_T=\es$. Using this, \er{mafla} and the identity $\Phi_{\wt\gm}=\Phi_{\mathbf{x}}$, we obtain
\[\lb{aeepe}
\wt\gm(\be)=\Phi_{\wt\gm}(\mathbf{c}_\be)=\Phi_{\mathbf{x}}(\mathbf{c}_\be),\qqq \forall\, \be\in\cS_T.
\]
Let $\be\in E_T$. Then $\wt\gm(\be)=0$. Thus, due to the definition \er{cat}, $\wt\gm=\gm(\,\cdot\,,T)$.

It remains to be shown that $T$ is a $\mathbf{x}$-minimal spanning tree of $G$. Let $\wt T$ be a spanning tree of $G$ such that $\b_{\wt T}<\b_T$. This and \er{sugm} yield that
$$
\textstyle\frac12\#\supp\gm(\,\cdot\,,\wt T)=\b_{\wt T}<\b_T=\frac12\#\supp\wt\gm,
$$
which contradicts the minimality of the form $\wt\gm$.

\emph{ii}) Let $T=(V,E_{T})$ be a $\mathbf{x}$-minimal spanning tree of the graph $G$. We prove that $\gm(\,\cdot\,,T)\in\mF(\mathbf{x})$ is a minimal form. The proof is by contradiction. Let $\wt\gm\in\mF(\mathbf{x})$ be a minimal form and
\[\lb{wtgm}
\#\supp\wt\gm<\#\supp\gm(\,\cdot\,,T).
\]
Then, due to item \emph{i}), there exists a $\mathbf{x}$-minimal spanning tree $\wt T$ of the graph $G$ such that
\[\lb{gmgb1}
\wt\gm=\gm(\,\cdot\,,\wt T).
\]
Combining \er{wtgm} and \er{gmgb1}, we obtain
\[\lb{wtTT}
\#\supp\gm(\,\cdot\,,\wt T)<\#\supp\gm(\,\cdot\,,T).
\]
This and the second identity in \er{sugm} yield $\b_{\wt T}<\b_T$,
where $\b_T$ is defined by \er{ncIT}, which contradicts the $\mathbf{x}$-minimality of the spanning tree $T$ and completes the proof.

\emph{iii}) Let $\wt\gm\in\mF(\mathbf{x})$ be a minimal form. Due to item \emph{i}) there exists a $\mathbf{x}$-minimal spanning tree $T$ of the graph $G$ such that
$\wt\gm=\gm(\,\cdot\,,T)$. Then we obtain using the identities \er{nbo} and \er{sugm}
\[\lb{smT}
\b(\mathbf{x})=\b_T=\textstyle\frac12\,\#\supp\gm(\,\cdot\,,T)=\frac12\,\#\supp\wt\gm.
\]

\setlength{\unitlength}{1.0mm}
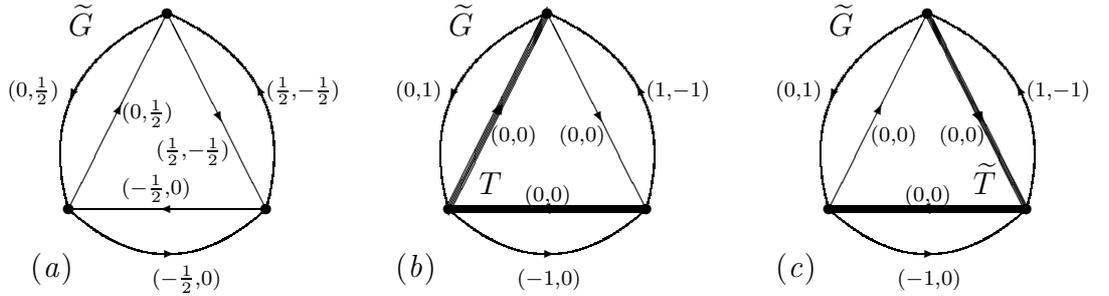
\begin{figure}[h]
\centering
\unitlength 1mm 
\linethickness{0.4pt}
\ifx\plotpoint\undefined\newsavebox{\plotpoint}\fi 
\hspace{-2.5cm}
\begin{picture}(146,45)(0,0)

\put(20,6){(\emph{a})}
\put(25,38){$\wt G$}
\put(31.0,27.0){\vector(1,2){1.0}}
\put(44.2,28.5){\vector(1,-2){1.0}}
\put(38,15.0){\vector(-1,0){1.0}}

\put(25,15){\line(1,0){26.0}}
\put(25,15){\line(1,2){13.0}}

\put(50.8,29.5){\vector(-1,2){1.0}}
\put(26.2,31.5){\vector(-1,-2){1.0}}
\put(38,9.0){\vector(1,0){1.0}}

\put(51,15){\line(-1,2){13.0}}

\put(25,15){\circle*{1.5}}
\put(51,15){\circle*{1.5}}
\put(38,41){\circle*{1.5}}
\bezier{200}(25,15)(38,3)(51,15)
\bezier{200}(25,15)(20,33)(38,41)
\bezier{200}(51,15)(56,33)(38,41)

\put(32.0,17.0){$\scriptstyle(-{1\/2},0)$}
\put(32.0,27.0){$\scriptstyle(0,{1\/2})$}
\put(36.5,22.0){$\scriptstyle({1\/2},-{1\/2})$}

\put(36.0,5.0){$\scriptstyle(-{1\/2},0)$}
\put(17.0,30.0){$\scriptstyle(0,{1\/2})$}
\put(51.0,30.0){$\scriptstyle({1\/2},-{1\/2})$}

\put(80.8,27.0){\vector(1,2){1.0}}
\put(81.0,27.0){\vector(1,2){1.0}}
\put(81.2,27.0){\vector(1,2){1.0}}
\put(94.2,28.5){\vector(1,-2){1.0}}
\put(88,15.2){\vector(-1,0){1.0}}
\put(88,15.0){\vector(-1,0){1.0}}
\put(88,14.8){\vector(-1,0){1.0}}

\put(75,14.6){\line(1,0){26.0}}
\put(75,14.8){\line(1,0){26.0}}
\put(75,15){\line(1,0){26.0}}
\put(75,15.2){\line(1,0){26.0}}
\put(75,15.4){\line(1,0){26.0}}

\put(75,14.2){\line(1,2){13.0}}
\put(75,14.6){\line(1,2){13.0}}
\put(75,15){\line(1,2){13.0}}
\put(75,15.4){\line(1,2){13.0}}
\put(75,15.8){\line(1,2){13.0}}

\put(100.8,29.5){\vector(-1,2){1.0}}
\put(76.2,31.5){\vector(-1,-2){1.0}}
\put(88,9.0){\vector(1,0){1.0}}

\put(101,15){\line(-1,2){13.0}}

\put(75,15){\circle*{1.5}}
\put(101,15){\circle*{1.5}}
\put(88,41){\circle*{1.5}}
\bezier{200}(75,15)(88,3)(101,15)
\bezier{200}(75,15)(70,33)(88,41)
\bezier{200}(101,15)(106,33)(88,41)

\put(85.0,16.0){$\scriptstyle(0,0)$}
\put(80.5,24.0){$\scriptstyle(0,0)$}
\put(89.5,24.0){$\scriptstyle(0,0)$}

\put(84.0,5.0){$\scriptstyle(-1,0)$}
\put(68.0,30.0){$\scriptstyle(0,1)$}
\put(101.0,30.0){$\scriptstyle(1,-1)$}

\put(79.0,17.0){$T$}
\put(68.0,6.0){(\emph{b})}
\put(75,38){$\wt G$}
\put(131.0,27.0){\vector(1,2){1.0}}
\put(144.4,28.5){\vector(1,-2){1.0}}
\put(144.2,28.5){\vector(1,-2){1.0}}
\put(144.0,28.5){\vector(1,-2){1.0}}
\put(138,15.2){\vector(-1,0){1.0}}
\put(138,15.0){\vector(-1,0){1.0}}
\put(138,14.8){\vector(-1,0){1.0}}

\put(125,14.6){\line(1,0){26.0}}
\put(125,14.8){\line(1,0){26.0}}
\put(125,15){\line(1,0){26.0}}
\put(125,15.2){\line(1,0){26.0}}
\put(125,15.4){\line(1,0){26.0}}

\put(125,15){\line(1,2){13.0}}

\put(150.8,29.5){\vector(-1,2){1.0}}
\put(126.2,31.5){\vector(-1,-2){1.0}}
\put(138,9.0){\vector(1,0){1.0}}

\put(151,14.4){\line(-1,2){13.0}}
\put(151,14.7){\line(-1,2){13.0}}
\put(151,15){\line(-1,2){13.0}}
\put(151,15.3){\line(-1,2){13.0}}
\put(151,15.6){\line(-1,2){13.0}}

\put(125,15){\circle*{1.5}}
\put(151,15){\circle*{1.5}}
\put(138,41){\circle*{1.5}}
\bezier{200}(125,15)(138,3)(151,15)
\bezier{200}(125,15)(120,33)(138,41)
\bezier{200}(151,15)(156,33)(138,41)

\put(118.0,6.0){(\emph{c})}
\put(125,38){$\wt G$}
\put(144.0,17.0){$\wt T$}
\put(135.0,16.0){$\scriptstyle(0,0)$}
\put(130.5,24.0){$\scriptstyle(0,0)$}
\put(139.5,24.0){$\scriptstyle(0,0)$}

\put(134.0,5.0){$\scriptstyle(-1,0)$}
\put(118.0,30.0){$\scriptstyle(0,1)$}
\put(151.0,30.0){$\scriptstyle(1,-1)$}
\end{picture}

\vspace{-0.7cm} \caption{\footnotesize  \emph{a}) A graph $\wt G$; the values of the 1-form $\wt{\mathbf{x}}$ on it are shown near edges; \; \emph{b}),\,\emph{c}) distinct $\wt{\mathbf{x}}$-minimal spanning trees $T$ and $\wt T$ (their edges are marked by bold) correspond to the same minimal form $\wt\gm\in\mF(\wt{\mathbf{x}})$, where $\wt\gm=\gm(\,\cdot\,,T)=\gm(\,\cdot\,,\wt T\,)$. The values of $\wt\gm$ are shown near edges.}
\label{TwtT}
\end{figure}

\emph{iv}) We consider the graph $\wt G$ shown in Fig.\ref{TwtT}\emph{a} and define the 1-form $\wt{\mathbf{x}}$ on it as shown in this figure. Looking over all spanning trees of $\wt G$ (the number of spanning trees is finite) one can check that the spanning trees $T$ (Fig.\ref{TwtT}\emph{b}) and $\wt T$ (Fig.\ref{TwtT}\emph{c}) are $\wt{\mathbf{x}}$-minimal spanning trees of the graph $\wt G$. By a direct calculation we obtain $\gm(\,\cdot\,,T)=\gm(\,\cdot\,,\wt T\,)$. \qq $\BBox$

\medskip

\no {\bf Remarks.} 1) Theorem \ref{Pphi}.\emph{i}) -- \emph{ii}) gives a full description of minimal forms $\gm\in\mF(\mathbf{x})$. Each minimal form $\gm$ is given by \er{cat} for some $\mathbf{x}$-minimal spanning tree $T$ and, conversely, for each $\mathbf{x}$-minimal spanning tree $T$ the form $\gm(\,\cdot\,,T)$ defined by \er{cat} is minimal.

2)  Due to Theorem \ref{Pphi}.\emph{iv}), the correspondence between minimal forms in a fixed set $\mF(\mathbf{x})$ on a graph $G$ and $\mathbf{x}$-minimal spanning trees of $G$ is not a bijection. Two different spanning trees may be associated with the same minimal form. Thus, the number of minimal forms in the set $\mF(\mathbf{x})$ is not greater than the number of $\mathbf{x}$-minimal spanning trees of the graph $G$.

\section{\lb{Sec3} Direct integrals}
\setcounter{equation}{0}

In this section we consider some properties of the coordinate form $\k$ on the fundamental graph $G_*$ defined by \er{edco}, \er{dco} and introduce an important 1-form $\t\in\mF(\k)$ on the graph $G_*$ called the \emph{index form}. We also give a general representation of fiber Laplacians in terms of any 1-form $\gb\in\mF(\k)$. Then we prove Theorem \ref{TDImf} about the decomposition for Schr\"odinger operators in terms of minimal forms and Proposition \ref{TNNI} about some properties of the invariant $\cI$ given by \er{dIm}.

\subsection{Properties of coordinate forms}
We formulate some properties of the coordinate form $\k$ on the fundamental graph $G_*$ defined by \er{edco}, \er{dco}.

\begin{proposition}
\lb{Pal0}
i) The coordinate form $\k:\cA_*\ra\R^d$ defined by \er{edco}, \er{dco} is a maximal form on the fundamental graph $G_*$ of the $\G$-periodic graph $\cG$ and the image of its flux function $\Phi_{\k}$ defined by \er{mafla} is the lattice $\Z^d$:
\[\lb{prao}
\Phi_{\k}(\cC)=\Z^d.
\]

ii) The kernel of the flux function $\Phi_{\k}$ does not depend on the choice of

$\bu$ the embedding of $\cG$ into $\R^d$;

$\bu$ the basis $a_1,\ldots,a_d$ of the lattice $\G$,

\no i.e., $\ker\Phi_{\k}$
is an invariant of the periodic graph $\cG$.

iii) The dimension of the kernel of $\Phi_{\k}$ satisfies
\[\lb{dike}
\dim\ker\Phi_{\k}=\b-d,
\]
where $\b$ is the Betti number of the fundamental graph $G_*$ given in \er{nni}.
\end{proposition}

\no \textbf{Remark.} Any cycle $\textbf{c}$ on the fundamental graph $G_*$ is obtained by factorization of a path on the periodic graph $\cG$ connecting some $\G$-equivalent vertices $v\in\cV$ and $v+a\in\cV$, $a\in\G$. Furthermore, the flux $\Phi_{\k}(\textbf{c})$ of the coordinate form $\k$ through the cycle $\mathbf{c}$ is equal to $\mn=(n_1,\ldots,n_d)\in\Z^d$, where $a=n_1a_1+\ldots+n_da_d$. In particular, $\Phi_{\k}(\textbf{c})=0$ if and only if the cycle $\textbf{c}$ on $G_*$ corresponds to a cycle on $\cG$.

\medskip

\no \textbf{Proof of Proposition \ref{Pal0}.}
\emph{i}) Obviously, $\k(\ul\be)=-\k(\be)$ for all $\be\in\cA_*$. Then the vector-valued function $\k:\cA_*\ra\R^d$ is a 1-form on $G_*$. All edges of the fundamental graph except loops corresponding to loops on the periodic graph have non-zero coordinates. Since $\Phi_{\k}=\Phi_{\gb}$ for any 1-form $\gb\in\mF(\k)$, due to the definition, $\k$ is a maximal form in $\mF(\k)$.

Now we prove \er{prao}. We take some vertex $v\in\cV$ of the periodic graph $\cG$.
Let $\mn=(n_1,\ldots,n_d)\in\Z^d$. Then $a=\sum\limits_{s=1}^dn_sa_s\in\G$. Since $\cG$ is connected, there exists an oriented path $\textbf{p}$ from $v$ to $v+a$ on $\cG$. Then $\textbf{c}=\textbf{p}/\G$ is a cycle on the fundamental graph $G_*=\cG/\G$ with the flux $\Phi_{\k}(\textbf{c})=\mn$. Thus, $\mn\in\Phi_{\k}(\cC)$. Conversely, let $\mn=(n_1,\ldots,n_d)\in\R^d$ and $\mn=\Phi_{\k}(\textbf{c})$ for some cycle $\textbf{c}\in\cC$. Then there exists an oriented path $\textbf{p}$ on the periodic graph $\cG$ such that $\textbf{c}=\textbf{p}/\G$ and $\textbf{p}$ connects the $\G$-equivalent vertices $v$ and $v+a$ on $\cG$, where $a=\sum\limits_{s=1}^dn_sa_s$. Thus, $a\in\G$ and, consequently, $\mn\in\Z^d$.

\emph{ii}) The kernel of the function $\Phi_{\k}$ consists of only cycles of the fundamental graph which correspond to cycles on the periodic graph $\cG$. But for any embedding of $\cG$ into $\R^d$ and for any basis of the lattice $\G$ any cycle on $\cG$ has zero sum of its edge coordinates. Thus, $\ker\Phi_{\k}$ does not depend on the graph embedding and the basis of $\G$.

\emph{iii}) From \er{prao} it follows that $\dim\Phi_{\k}(\cC)=d$. Then, using \er{debn}, we obtain
$$
\dim\ker\Phi_{\k}=\dim\cC-\dim\Phi_{\k}(\cC)=\b-d.
$$
\qq \BBox

\subsection{Edge indices.} In order to prove Theorem \ref{TDImf} about the direct integral decomposition we need to define an {\it edge index}, which was introduced in \cite{KS14}. The indices are important to study the spectrum of the Laplacians and Schr\"odinger operators on periodic graphs, since the initial fiber operators are expressed in terms of indices of the fundamental graph edges (see \er{l2.15''}).

For any vertex $v\in\cV$ of a $\G$-periodic graph $\cG$ the following unique representation holds true:
\[
\lb{Dv} v=\{v\}+[v], \qq \textrm{where}\qq \{v\}\in \Omega,\qquad [v]\in\G,
\]
$\Omega$ is the fundamental cell of the lattice $\G$ defined by \er{fuce}.
In other words, each vertex $v$ can be obtained from a vertex $\{v\}\in \Omega$ by a shift by the vector $[v]\in\G$. We call $\{v\}$ and $[v]$ the \emph{fractional and integer parts of the vertex $v$}, respectively.
For any oriented edge $\be=(u,v)\in\cA$ we define the {\bf edge "index"}
$\t(\be)$ as the vector of the lattice $\Z^d$ given by
\[
\lb{in}
\t(\be)=[v]_\A-[u]_\A\in\Z^d,
\]
where $[v]\in\G$ is the integer part of the vertex $v$ and the vector $[v]_\A\in\Z^d$ is defined by \er{cola}.

For example, for the periodic graph $\cG$ shown in Fig.\ref{ff.0.11}\emph{a} the index of the edge $(v_3+a_2,v_2+a_1)$ is equal to $(1,-1)$, since the integer parts of the vertices $v_3+a_2$ and $v_2+a_1$ are equal to $a_2$ and $a_1$, respectively.

\setlength{\unitlength}{1.0mm}
\begin{figure}[h]
\centering
\unitlength 1mm 
\linethickness{0.4pt}
\ifx\plotpoint\undefined\newsavebox{\plotpoint}\fi 

\begin{picture}(160,55)(0,0)
\bezier{30}(81,15)(87.5,28)(94,41)
\bezier{30}(82,15)(88.5,28)(95,41)
\bezier{30}(83,15)(89.5,28)(96,41)
\bezier{30}(84,15)(90.5,28)(97,41)
\bezier{30}(85,15)(91.5,28)(98,41)
\bezier{30}(86,15)(92.5,28)(99,41)
\bezier{30}(87,15)(93.5,28)(100,41)
\bezier{30}(88,15)(94.5,28)(101,41)
\bezier{30}(89,15)(95.5,28)(102,41)
\bezier{30}(90,15)(96.5,28)(103,41)
\bezier{30}(91,15)(97.5,28)(104,41)
\bezier{30}(92,15)(98.5,28)(105,41)
\bezier{30}(93,15)(99.5,28)(106,41)
\bezier{30}(94,15)(100.5,28)(107,41)
\bezier{30}(95,15)(101.5,28)(108,41)
\bezier{30}(96,15)(102.5,28)(109,41)
\bezier{30}(97,15)(103.5,28)(110,41)
\bezier{30}(98,15)(104.5,28)(111,41)
\bezier{30}(99,15)(105.5,28)(112,41)
\bezier{30}(100,15)(106.5,28)(113,41)
\bezier{30}(101,15)(107.5,28)(114,41)
\bezier{30}(102,15)(108.5,28)(115,41)
\bezier{30}(103,15)(109.5,28)(116,41)
\bezier{30}(104,15)(110.5,28)(117,41)
\bezier{30}(105,15)(111.5,28)(118,41)

\bezier{25}(21,15)(26,25)(31,35)
\bezier{25}(22,15)(27,25)(32,35)
\bezier{25}(23,15)(28,25)(33,35)
\bezier{25}(24,15)(29,25)(34,35)
\bezier{25}(25,15)(30,25)(35,35)
\bezier{25}(26,15)(31,25)(36,35)
\bezier{25}(27,15)(32,25)(37,35)
\bezier{25}(28,15)(33,25)(38,35)
\bezier{25}(29,15)(34,25)(39,35)
\bezier{25}(30,15)(35,25)(40,35)
\bezier{25}(31,15)(36,25)(41,35)
\bezier{25}(32,15)(37,25)(42,35)
\bezier{25}(33,15)(38,25)(43,35)
\bezier{25}(34,15)(39,25)(44,35)
\bezier{25}(35,15)(40,25)(45,35)
\bezier{25}(36,15)(41,25)(46,35)
\bezier{25}(37,15)(42,25)(47,35)
\bezier{25}(38,15)(43,25)(48,35)
\bezier{25}(39,15)(44,25)(49,35)

\put(33.0,25){$\Omega$}

\put(95,25){$\Omega$}

\put(0,15){\line(1,0){60.0}}
\put(10,35){\line(1,0){60.0}}
\put(20,55){\line(1,0){60.0}}
\put(0,15){\line(1,2){20.0}}
\put(20,15){\line(1,2){20.0}}
\put(40,15){\line(1,2){20.0}}
\put(60,15){\line(1,2){20.0}}

\put(10,15){\line(-1,2){5.0}}
\put(30,15){\line(-1,2){15.0}}
\put(50,15){\line(-1,2){20.0}}
\put(65,25){\line(-1,2){15.0}}
\put(75,45){\line(-1,2){5.0}}

\put(20,15){\vector(1,0){20.0}}
\put(20,15){\vector(1,2){10.0}}

\put(0,15){\circle{1}}
\put(10,15){\circle{1}}
\put(20,15){\circle*{1.5}}
\put(30,15){\circle*{1.5}}
\put(40,15){\circle{1}}
\put(50,15){\circle{1}}
\put(60,15){\circle{1}}

\put(5,25){\circle{1}}
\put(25,25){\circle*{1.5}}
\put(45,25){\circle{1}}
\put(65,25){\circle{1}}

\put(10,35){\circle{1}}
\put(20,35){\circle{1}}
\put(30,35){\circle{1}}
\put(40,35){\circle{1}}
\put(50,35){\circle{1}}
\put(60,35){\circle{1}}
\put(70,35){\circle{1}}

\put(15,45){\circle{1}}
\put(35,45){\circle{1}}
\put(55,45){\circle{1}}
\put(75,45){\circle{1}}

\put(20,55){\circle{1}}
\put(30,55){\circle{1}}
\put(40,55){\circle{1}}
\put(50,55){\circle{1}}
\put(60,55){\circle{1}}
\put(70,55){\circle{1}}
\put(80,55){\circle{1}}

\put(30.7,11.5){$\scriptstyle a_1$}
\put(23.5,31.0){$\scriptstyle a_2$}
\put(18.0,11.5){$\scriptstyle v_1$}
\put(20.0,37.0){$\scriptstyle v_1+a_2$}
\put(39.0,36.5){$\scriptstyle v_3+a_2$}
\put(38.0,11.5){$\scriptstyle v_1+a_1$}
\put(20.0,25.5){$\scriptstyle v_2$}
\put(46.0,24.5){$\scriptstyle v_2+a_1$}
\put(49.0,32.0){$\scriptstyle v_1+a_1+a_2$}
\put(29.5,16.5){$\scriptstyle v_3$}
\put(-10,15.0){(\emph{a})}


\put(80,15){\vector(1,0){26.0}}
\put(80,15){\vector(1,2){13.0}}
\bezier{30}(93,41)(106,41)(119,41)
\bezier{30}(106,15)(112.5,28)(119,41)
\put(80,15){\circle*{1.5}}
\put(93,15){\circle*{1.5}}
\put(106,15){\circle{1}}

\put(93,15){\line(-1,2){6.5}}
\put(112.5,28){\line(-1,2){6.5}}

\put(86.5,28){\circle*{1.5}}
\put(112.5,28){\circle{1}}

\put(93,41){\circle{1}}
\put(106,41){\circle{1}}
\put(119,41){\circle{1}}

\put(96.7,12.0){$a_1$}
\put(84.5,34.5){$a_2$}
\put(78.0,11.0){$v_1$}
\put(91.0,42.7){$v_1$}
\put(102.0,42.7){$v_3$}
\put(104.0,11.0){$v_1$}
\put(80.0,28.0){$v_2$}
\put(106.0,28.0){$v_2$}
\put(116.0,42.7){$v_1$}
\put(90.0,11.0){$v_3$}
\put(70,15.0){(\emph{b})}

\put(131.0,27.0){\vector(1,2){1.0}}
\put(144.2,28.5){\vector(1,-2){1.0}}
\put(138,15.0){\vector(-1,0){1.0}}

\put(125,15){\line(1,0){26.0}}

\put(125,15){\line(1,2){13.0}}

\put(150.8,29.5){\vector(-1,2){1.0}}
\put(126.2,31.5){\vector(-1,-2){1.0}}
\put(138,9.0){\vector(1,0){1.0}}

\put(151,15){\line(-1,2){13.0}}

\put(125,15){\circle*{1.5}}
\put(119.0,13.5){$v_1$}
\put(151,15){\circle*{1.5}}
\put(151.5,13.5){$v_3$}
\put(138,41){\circle*{1.5}}
\put(138.0,42){$v_2$}
\bezier{200}(125,15)(138,3)(151,15)
\bezier{200}(125,15)(120,33)(138,41)
\bezier{200}(151,15)(156,33)(138,41)

\put(135.0,16.0){$\scriptstyle(0,0)$}
\put(130.5,24.0){$\scriptstyle(0,0)$}
\put(139.5,24.0){$\scriptstyle(0,0)$}

\put(134.0,5.0){$\scriptstyle(-1,0)$}
\put(123.0,37.0){$\scriptstyle(0,1)$}
\put(145.0,37.0){$\scriptstyle(1,-1)$}
\end{picture}

\vspace{-0.7cm} \caption{\footnotesize \emph{a}) The Kagome lattice, the fundamental cell $\Omega$ is shaded; \emph{b})~ its fundamental graph with a minimal form $\gm\in\mF(\k)$. The values of $\gm$ are shown near edges.}
\label{ff.0.11}
\end{figure}
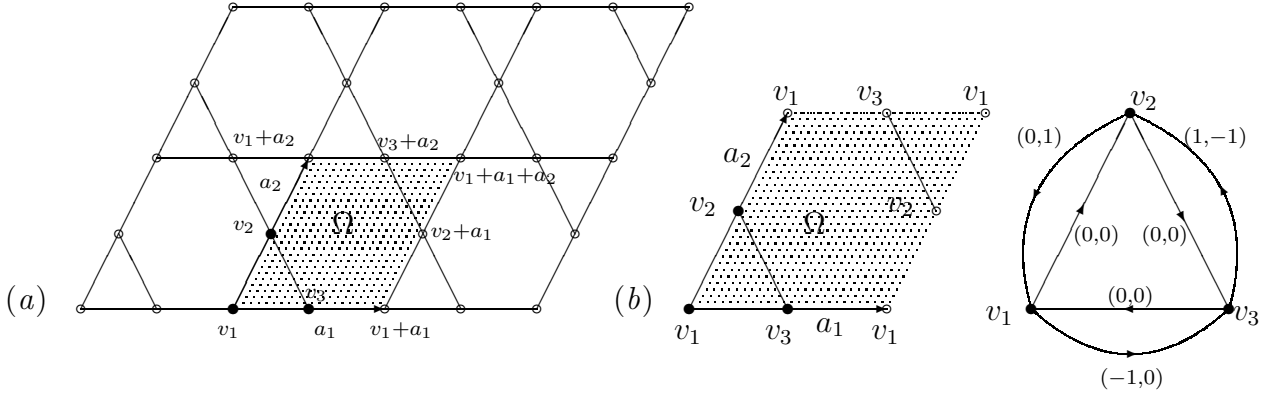

On the fundamental graph $G_*$ we introduce the \emph{index form} $\t: \cA_*\ra\Z^d$ by:
\[
\lb{inf}
\t(\bf e_*)=\t(\be) \qq \textrm{ for some $\be\in\cA$ \; such that }  \; \be_*=\gf(\be), \qqq \be_*\in\cA_*,
\]
where $\gf$ is defined by \er{sur}.

The index form $\t$ is uniquely determined by \er{inf}, since
$$
\t(\be+a)=\t(\be),\qqq \forall\, (\be,a)\in\cA \ts\G.
$$

\begin{proposition}
\lb{Ptau}
The index form $\tau:\cA_*\ra\Z^d$ defined by \er{in}, \er{inf} is a 1-form on the fundamental graph $G_*$ and
\[\lb{PaPt}
\Phi_{\tau}=\Phi_{\k},
\]
where $\Phi_{\gb}$ is the flux function of a 1-form $\gb$ defined by \er{mafla}, $\k$ is the coordinate form given by \er{edco}, \er{dco}.
\end{proposition}

\no \textbf{Proof.} From the definitions \er{in}, \er{inf} of the
index form we have
\[\lb{inin}
\t(\ul\be)=-\t(\be) \qq \textrm{for each }\qq\be\in\cA_*.
\]
Thus, the function $\t:\cA_*\ra\Z^d$ is a 1-form on $G_*$. Now we show \er{PaPt}.
Let $\textbf{c}\in\cC$, where $\cC$ is the cycle space of $G_*$. Then there exists an oriented path $\textbf{p}=(\be_1,\be_2,\ldots,\be_n)$, $\be_s=(v_s,v_{s+1})$, $s\in\N_n$, on the periodic graph $\cG$ such that $\textbf{c}=\textbf{p}/\G$ and $\textbf{p}$ connects the $\G$-equivalent vertices $v_1$ and $v_{n+1}=v_1+a$ on $\cG$ for some $a\in\G$. Then, due to the definition \er{edco}, \er{dco} of the coordinate form $\k$, we have
$$
\begin{aligned}
&\Phi_{\k}(\textbf{c})=\sum_{s=1}^n\k(\be_s)=(v_2)_\A-(v_1)_\A+(v_3)_\A-(v_2)_\A+
\ldots+(v_{n+1})_\A-(v_n)_\A\\
&=(v_{n+1})_\A-(v_1)_\A=(v_1+a)_\A-(v_1)_\A=a_\A,
\end{aligned}
$$
where $x_\A$ is defined by \er{cola}. On the other hand, due to the definition \er{in}, \er{inf} of the index form, we obtain
$$
\begin{aligned}
&\Phi_{\t}(\textbf{c})=\sum_{s=1}^n\t(\be_s)=[v_2]_\A-[v_1]_\A+[v_3]_\A-[v_2]_\A+
\ldots+[v_{n+1}]_\A-[v_n]_\A\\
&=[v_{n+1}]_\A-[v_1]_\A=[v_1+a]_\A-[v_1]_\A=a_\A,
\end{aligned}
$$
where $[v]$ is the integer part of the vertex $v$. Thus, $\Phi_{\k}(\textbf{c})=\Phi_{\t}(\textbf{c})$. Since $\textbf{c}$ is any cycle on $G_*$, \er{PaPt} has been proved. \qq \BBox

\medskip

\no \textbf{Remark.}  For some periodic graphs the index form $\t\in\mF(\k)$ may coincide with maximal or (and) minimal forms in $\mF(\k)$.

\subsection{Direct integral decomposition.} Recall that we introduce the Hilbert space $\mH$ by \er{Hisp}, $\lan\cdot\,,\cdot\ran$ denotes the standard inner product in $\R^d$, and $\{a_1,\ldots,a_d\}$ is the basis of the lattice $\G$. We identify the vertices of the fundamental
graph $G_*=(\cV_*,\cE_*)$ with the vertices of the $\G$-periodic graph
$\cG=(\cV,\cE)$ from the fundamental cell $\Omega$. We need the following results (see Theorem 1.1 in \cite{KS14}).

\begin{theorem}\label{TFD1}
The Schr\"odinger operator $H=\D+Q$ on $\ell^2(\cV)$ has the following decomposition into a constant fiber direct integral
\[
\lb{raz}
\begin{aligned}
& UH U^{-1}=\int^\oplus_{\T^d}H_\t(\vt){d\vt\/(2\pi)^d}
\end{aligned}
\]
for the unitary Gelfand transform $U:\ell^2(\cV)\to\mH$ defined by
\[
\lb{5001}
(Uf)(\vt,v)=\sum\limits_{\mn=(n_1,\ldots,n_d)\in\Z^d}e^{-i\lan \mn,\vt\ran }
f(v+n_1a_1+\ldots+n_da_d), \qqq (\vt,v)\in \T^d\ts\cV_*.
\]
Here the fiber Schr\"odinger operator $H_\t(\vt)$ and the fiber Laplacian $\D_\t(\vt)$ are given by
\[\label{Hvt'}
H_\t(\vt)=\D_\t(\vt)+Q, \qqq \forall\,\vt\in\T^d,
\]
\[
\label{l2.15''}
 \big(\D_\t(\vt)f\big)(v)=\vk_vf(v)-\sum_{\be=(v,\,u)\in\cA_*} e^{i\lan\t(\be),\,\vt\ran}f(u), \qqq f\in\ell^2(\cV_*),\qqq
 v\in\cV_*,
\]
where $\t$ is the index form defined by \er{in}, \er{inf}.
\end{theorem}

The fiber Laplacians $\D_\t(\cdot)$ given by \er{l2.15''} depend on $\#\cA_*$ indices $\t(\be)$ of the fundamental graph edges $\be\in\cA_*$. Some of them may be zero, but the number of such zero indices is difficult to control and depends on the choice of the embedding of the periodic graph $\cG$ into the space $\R^d$. Therefore, in Theorem \ref{TDIf} we construct some gauge transformation and give a representation of the fiber Laplacians in terms of any 1-form $\gb\in\mF(\k)$. In particular, if $\gb$ is a minimal form, then the fiber operator has the minimal number of coefficients depending on the quasimomentum and it will be used to estimate the Lebesgue measure of the spectrum of the Schr\"odinger operators on periodic graphs (see Theorem \ref{T1}) and the effective masses at the bottom of the Laplacian spectrum (see Theorem \ref{TTsem}).

\medskip

Let $\gb\in\mF(\k)$ be a 1-form on the fundamental graph $G_*$, where $\mF(\k)$ is given by \er{vv1f} at $\mathbf{x}=\k$, and $\k$ is the coordinate form defined by \er{edco}, \er{dco}. We fix a vertex $v_0\in\cV_*$ and introduce the gauge transformation $W_\gb:\mH\to\mH$ by
\[\lb{Uvt}
(W_\gb \,g)(\vt,v)=e^{i\lan w_\gb(v),\vt\ran}g(\vt,v), \qqq g\in\mH,\qqq  (\vt,v)\in \T^d\ts\cV_*,
\]
where the function $w_\gb:\cV_*\ra\R^d$ is defined as follows: for any vertex
$v\in\cV_*$, we take an oriented path $\textbf{p}=(\be_1,\be_2,\ldots,\be_n)$ on $G_*$
starting at $v_0$ and ending at $v$ and  we set
\[\lb{Wvt}
w_\gb(v)=\sum_{s=1}^n\big(\t(\be_s)-\gb(\be_s)\big).
\]

We show that $w_\gb(v)$ does not depend on the choice of a path from $v_0$ to $v$. Indeed, let $\textbf{p}$ and $\textbf{q}$ be some oriented paths from
$v_0$ to $v$. We consider the cycle $\mathbf{c}=\textbf{p}\ul{\textbf{q}}$, where $\ul{\textbf{q}}$ is the inverse path of $\textbf{q}$. Then we have
\[\lb{com1}
\sum_{\be\in\mathbf{c}}\big(\t(\be)-\gb(\be)\big)=\sum_{\be\in
\textbf{p}}\big(\t(\be)-\gb(\be)\big)-\sum_{\be\in
\textbf{q}}\big(\t(\be)-\gb(\be)\big).
\]
From the identities $\Phi_\gb=\Phi_{\k}$ and \er{PaPt} it follows that
\[\lb{PbPt}
\Phi_{\gb}=\Phi_{\t}.
\]
Combining the identities \er{com1} and \er{PbPt}, we obtain
$$
\sum_{\be\in \mathbf{p}}\big(\t(\be)-\gb(\be)\big)=\sum_{\be\in
\mathbf{q}}\big(\t(\be)-\gb(\be)\big),
$$
which implies that $w_\gb(v)$ does not depend on the choice of the path from $v_0$ to $v$.
\begin{theorem}
\label{TDIf}
i) Let $\gb\in\mF(\k)$, where $\k$ is the coordinate form defined by \er{edco}, \er{dco}. Then the Schr\"odinger operator $H=\D+Q$ on $\ell^2(\cV)$ has the following decomposition into a constant fiber direct integral
\[
\lb{raz1}
\mU_\gb H \mU_\gb^{-1}=\int^\oplus_{\T^d}H_\gb(\vt)\frac{d\vt
}{(2\pi)^d}\,,
\]
where the unitary operator $\mU_\gb=W_\gb U:\ell^2(\cV)\to\mH$, $U$ is the Gelfand transformation given by \er{5001} and $W_\gb$ is the gauge transformation
defined by \er{Uvt}, \er{Wvt}.  Here the fiber Schr\"odinger operator $H_\gb(\vt)$ and the fiber Laplacian $\D_\gb(\vt)$ are given by
\[
\label{Hvt}
H_\gb(\vt)=\D_\gb(\vt)+Q,\qqq \forall\,\vt\in\T^d,
\]
\[
\label{l2.13a}
\begin{aligned}
\big(\D_\gb(\vt)f\big)(v)=\vk_vf(v)
-\sum_{\be=(v,u)\in\cA_*}e^{i\lan\gb(\be),\,\vt\ran}f(u), \qqq f\in\ell^2(\cV_*), \qqq v\in\cV_*.
\end{aligned}
\]

ii) Let $\gb'$ be another form from the set $\mF(\k)$. Then for each $\vt\in\T^d$ the operators $\D_\gb(\vt)$ and $\D_{\gb'}(\vt)$ are unitarily equivalent.

iii) The number $\cN_\gb$ of exponents $e^{i\lan\gb(\be),\,\cdot\,\ran}\neq1$, $\be\in\cA_*$, in the identities \er{l2.13a} for the fiber Laplacians $\D_\gb(\cdot)$ satisfies
\[\lb{dI}
\cN_\gb\geq2\cI, \qqq \textrm{where} \qqq \cI=\textstyle\frac12\,\#\supp \gm=\b(\k),
\]
$\gm\in\mF(\k)$ is a minimal form, and $\b(\k)$ is defined by \er{nbo} at $\mathbf{x}=\k$.
\end{theorem}

\no{\bf Proof.} \emph{i}) We show that the operators $\D_\t(\cdot)$ and $\D_\gb(\cdot)$ defined by \er{l2.15''} and \er{l2.13a}, respectively, are unitarily equivalent, by the gauge transformation $W_\gb$ given by \er{Uvt}. Indeed, from \er{Wvt} it follows that
\[\lb{waa}
w_\gb(u)=w_\gb(v)+\t(\be)-\gb(\be), \qqq \forall\,\be=(v,u)\in\cA_*.
\]
Using this, \er{l2.15''}, \er{Uvt} and \er{l2.13a}, we have
\[
\begin{aligned}
&\big(\D_\t(\vt)f\big)(v)=\vk_vf(v)-\sum_{\be=(v,u)\in\cA_*}
e^{i\lan\t(\be),\vt\ran}f(u)\\
&=\vk_vf(v)-\sum_{\be=(v,u)\in\cA_*}e^{-i\lan w_\gb(v),\vt\ran}e^{i\lan\gb(\be),\vt\ran}e^{i\lan w_\gb(u),\vt\ran}f(u)\\
&=e^{-i\lan w_\gb(v),\vt\ran}\Big(\vk_ve^{i\lan w_\gb(v),\vt\ran}f(v)
-\sum_{\be=(v,u)\in\cA_*}e^{i\lan\gb(\be),\vt\ran}e^{i\lan w_\gb(u),\vt\ran}f(u)\Big)\\
&=\big(W_\gb^{-1}\D_{\gb}(\vt)W_\gb f\big)(\vt,v).
\end{aligned}
\]
Thus,
\[\lb{DaDa}
\D_\t(\cdot)=W_\gb^{-1}\D_{\gb}(\cdot)W_\gb,
\]
i.e., $\D_\gb(\vt)$ and $\D_\t(\vt)$ are unitarily equivalent for each $\vt\in\T^{d}$. Consequently, the fiber Schr\"odinger operator $H_\t(\cdot)=\D_\t(\cdot)+Q$ is
unitarily equivalent to the operator $H_\gb(\cdot)=\D_\gb(\cdot)+Q$, by the gauge transformation $W_\gb$. This and Theorem \ref{TFD1} give the required statement.

\emph{ii}) This is a direct consequence of the identity \er{DaDa}.

\emph{iii}) Due to the definition of the number $\cN_\gb$ we have $\cN_\gb=\#\supp \gb$ for each 1-form $\gb\in\mF(\k)$. Since $\gm\in\mF(\k)$ is a minimal form,
\[\lb{cNgb}
\cN_\gb=\#\supp \gb\geq \#\supp \gm.
\]
This, \er{dIm} and Theorem \ref{Pphi}.\emph{iii}) give \er{dI}. \qq \BBox

\medskip

\no \textbf{Remarks.} 1) The formulas \er{raz1} -- \er{l2.13a} give an infinite number of decompositions into constant fiber direct integrals for the same Schr\"odinger operator $H$ on the periodic graph $\cG$ (one decomposition for each $\gb\in\mF(\k)$).

2) From Theorem \ref{TDIf}.\emph{iii}) it follows that the fiber Laplacians $\D_\gm(\vt)$ with a minimal form $\gm\in\mF(\k)$ defined by \er{l2.13am} have the minimal number $2\cI$ of coefficients depending on the quasimomentum $\vt$ among all fiber Laplacians $\D_\gb(\vt)$, where $\gb\in\mF(\k)$.

\medskip

\no {\bf Proof of Theorem \ref{TDImf}.} \emph{i}) In virtue of Theorem \ref{Pphi}.\emph{i}), the minimal form $\gm\in\mF(\kappa)$ satisfies $\gm=\wt\gm(\,\cdot\,,T)$ for some $\k$-minimal spanning tree $T=(\cV_*,\cE_T)$ of the fundamental graph $G_*$, where $\wt\gm(\,\cdot\,,T):\cA_*\ra\R^d$ is defined by \er{cat} as $\mathbf{x}=\k$. Then for each $\be\in\cS_T=\cE_*\sm\cE_{T}$ we have
\[\lb{gmbe}
\gm(\be)=\wt\gm(\be,T)=\Phi_\k(\mathbf{c}_\be),
\]
where $\mathbf{c}_\be$ is a unique cycle on $G_*$ whose edges are all in $T$ except $\be$. Since the set of all cycles $\{\mathbf{c}_\be\}_{\be\in\cS_T}$ forms a basis of the cycle space $\cC$ of the graph $G_*$ and $\Phi_\k(\cC)=\Z^d$ (see \er{prao}), we conclude that $\gm(\cA_*)$ generates the group $\Z^d$. Consequently, there exist $d$ edges $\be_1,\ldots,\be_d\in\supp\gm$ such that the set
$$
\big\{\gm(\be_s)=\big(\gm_1(\be_s),\ldots,\gm_d(\be_s)\big)\in\Z^d: s\in\N_d\big\}
$$
forms a basis of $\Z^d$. Then the set
$$
\A'=\big\{a'_s=\gm_1(\be_s)a_1+\ldots+\gm_d(\be_s)a_d: s\in\N_d\big\}
$$
is a basis of the lattice $\G$. Thus, under the change of the basis of the lattice $\G$ from $\A$ to $\A'$ the values of the minimal form $\big\{\gm(\be_s)\big\}_{s\in\N_d}$ transform to an orthonormal basis of $\Z^d$.

\emph{ii}) Since $\gm\in\mF(\k)$, the decomposition \er{raz1m} -- \er{l2.13am} is a direct consequence of the decomposition \er{raz1} -- \er{l2.13a} as $\gb=\gm$.

\emph{iii}) Since $\gm,\gm'\in\mF(\k)$, this item follows from Theorem \ref{TDIf}.\emph{ii}) as $\gb=\gm$ and $\gb'=\gm'$.

\emph{iv}) This item is obvious.

\emph{v}) Due to Theorem \ref{Pphi}.\emph{iii}) we have
$\cI=\b_T$, where $\b_T$ is defined by \er{ncIT} as $\mathbf{x}=\k$, and $T$ is any $\k$-minimal spanning tree of the fundamental graph $G_*$. In virtue of Proposition \ref{Pal0}.\emph{ii}), $\ker\Phi_{\k}$ depends neither on the choice of the embedding of $\cG$ into the space $\R^d$ nor on
the choice of the basis of the lattice $\G$. Then the definition \er{ncIT} of the number $\b_T$ shows that $\cI$ does not depend on these choices. Due to the definition of the minimal form $\gm$, the number $\cI$ does not depend on the choice of $\gm\in\mF(\k)$.  \qq $\BBox$

\begin{corollary}\label{TCo2}
Let $\gm\in\mF(\k)$ be a minimal form on the fundamental graph $G_*=(\cV_*,\cE_*)$, where $\k$ is the coordinate form defined by \er{edco}, \er{dco}. Then in the standard orthonormal basis of $\ell^2(\cV_*)=\C^\n$, $\n=\#\cV_*$, the $\n\ts\n$ matrix $\D_\gm(\vt)=\big(\D_{\gm,uv}(\vt)\big)_{u,v\in\cV_*}$ of the fiber Laplacian $\D_\gm(\vt)$ defined by \er{l2.13am} is given by
\[\lb{Devt}
\D_{\gm,uv}(\vt)=\left\{
\begin{array}{ll}
\displaystyle\vk_v-\sum\limits_{\be=(u,u)\in\cA_*}\cos\lan\gm(\be),\vt\ran, & \textrm{ if } \, u=v\\[20pt]
   \displaystyle-\sum\limits_{\be=(u,v)\in\cA_*}
  e^{-i\lan\gm(\be),\vt\ran},  & \textrm{ if } \, u\sim v,\qq u\neq v\\[20pt]
 \hspace{20mm} 0, \qq & \textrm{ otherwise}\\
\end{array}\right..
\]
Here ${\vk}_v$ is the degree of the vertex $v$.
\end{corollary}

\no {\bf Proof.} Let $\{\gh_u\}_{u\in\cV_*}$ be the standard orthonormal basis of $\ell^2(\cV_*)$. Substituting the formula \er{l2.13am} into the identity
$$
\D_{\gm,uv}(\vt)=\lan \gh_u,\D_\gm(\vt)\gh_v\ran_{\ell^2(\cV_*)}
$$
and using the fact that for each loop $\be=(u,u)\in\cA_*$ there exists a loop $\ul\be=(u,u)\in\cA_*$ and $\gm(\ul\be)=-\gm(\be)$,
we obtain \er{Devt}. \qq $\BBox$

\medskip

\no \textbf{Remark.} The number of entries lying on and above the main diagonal of the matrix $\D_\gm(\vt)$, defined by \er{Devt}, and depending on $\vt$ can be strictly less than $\cI$, where $\cI=\textstyle\frac12\,\#\supp \gm$, since some edges from $\cE_*\cap\supp\gm$ may be incident to the same pair of vertices.

\medskip

\no \textbf{Proof of Proposition \ref{TNNI}.} Due to Theorem \ref{Pphi}.\emph{iii}) we have $\cI=\b_T$, where $\b_T$ is defined by \er{ncIT} as $\mathbf{x}=\k$, and $T$ is any $\k$-minimal spanning tree of the fundamental graph $G_*$. Using the definition \er{ncIT} of the number $\b_T$ and the identities \er{debn} and \er{dike}:
$$
\dim\cC=\b,\qqq \dim\ker\Phi_{\k}=\b-d,
$$
we obtain $d\leq \b_T\leq\b$. Thus, the estimate \er{nni} has been proved.

\emph{i}) -- \emph{ii}), \emph{v}) Let $\mn$ be any nonnegative integer number. We consider a finite connected graph $G_1$ with the Betti number $\b_1=\mn$. For the periodic graph $\cG$ obtained from the $d$-dimensional lattice by "gluing"\, the graph $G_1$ to each vertex of the lattice (see Fig. \ref{ff.10}), due to Proposition \ref{TG1}.\emph{i}), we have $d=\cI$ and the Betti number $\b$ of the fundamental graph $G_*$ satisfies
\[\lb{bdb0}
\b=d+\b_1=\cI+\mn.
\]
If $G_1$ is a tree, then $\b_1=0$. Consequently, $\b=d$ and we have the case \emph{i}). If $G_1$ is not a tree, then $\b_1>0$ and, consequently, $\b>d$ and we have the case \emph{ii}). Moreover, from \er{bdb0} it follows that $\b-\cI=\mn$.

\emph{iii}) For the triangular lattice (see Fig. \ref{ff.10T}) $d=2$, $\cI=\b=3$.

\emph{iv}) For the Kagome lattice (see Fig. \ref{ff.0.11}) $d=2$, $\cI=3$, $\b=4$. \qq \BBox

\setlength{\unitlength}{1.0mm}
\begin{figure}[h]
\centering
\unitlength 1.0mm 
\linethickness{0.4pt}
\ifx\plotpoint\undefined\newsavebox{\plotpoint}\fi 
\begin{picture}(120,45)(0,0)
\unitlength 1.0mm

\bezier{20}(20.5,20)(20.5,25)(20.5,30)
\bezier{20}(21.0,20)(21.0,25)(21.0,30)
\bezier{20}(21.5,20)(21.5,25)(21.5,30)
\bezier{20}(22.0,20)(22.0,25)(22.0,30)
\bezier{20}(22.5,20)(22.5,25)(22.5,30)
\bezier{20}(23.0,20)(23.0,25)(23.0,30)
\bezier{20}(23.5,20)(23.5,25)(23.5,30)
\bezier{20}(24.0,20)(24.0,25)(24.0,30)
\bezier{20}(24.5,20)(24.5,25)(24.5,30)
\bezier{20}(25.0,20)(25.0,25)(25.0,30)
\bezier{20}(25.5,20)(25.5,25)(25.5,30)
\bezier{20}(26.0,20)(26.0,25)(26.0,30)
\bezier{20}(26.5,20)(26.5,25)(26.5,30)
\bezier{20}(27.0,20)(27.0,25)(27.0,30)
\bezier{20}(27.5,20)(27.5,25)(27.5,30)
\bezier{20}(28.0,20)(28.0,25)(28.0,30)
\bezier{20}(28.5,20)(28.5,25)(28.5,30)
\bezier{20}(29.0,20)(29.0,25)(29.0,30)
\bezier{20}(29.5,20)(29.5,25)(29.5,30)
\bezier{20}(30.0,20)(30.0,25)(30.0,30)

\put(10,10){\line(1,0){30.00}}
\put(10,20){\line(1,0){30.00}}
\put(10,30){\line(1,0){30.00}}
\put(10,40){\line(1,0){30.00}}
\put(10,10){\line(0,1){30.00}}
\put(20,10){\line(0,1){30.00}}
\put(30,10){\line(0,1){30.00}}
\put(40,10){\line(0,1){30.00}}

\put(10,10){\line(1,1){30.00}}
\put(10,20){\line(1,1){20.00}}
\put(20,10){\line(1,1){20.00}}
\put(10,30){\line(1,1){10.00}}
\put(30,10){\line(1,1){10.00}}

\put(10,10){\circle{1}}
\put(20,10){\circle{1}}
\put(30,10){\circle{1}}
\put(40,10){\circle{1}}

\put(10,20){\circle{1.0}}
\put(20,20){\circle*{1.0}}
\put(30,20){\circle{1.0}}
\put(40,20){\circle{1}}

\put(20.5,18){$\scriptstyle v$}
\put(24.5,18){$\scriptstyle a_1$}
\put(16.5,25){$\scriptstyle a_2$}
\put(26.5,22){$\scriptstyle \Omega$}
\put(20,20){\vector(0,1){10.00}}
\put(20,20){\vector(1,0){10.00}}
\put(10,30){\circle{1}}
\put(20,30){\circle{1.0}}
\put(30,30){\circle{1.0}}
\put(40,30){\circle{1}}

\put(10,40){\circle{1}}
\put(20,40){\circle{1}}
\put(30,40){\circle{1}}
\put(40,40){\circle{1}}
\put(20,20){\line(1,1){10.00}}
\put(0,10){(\emph{a})}

\put(46,10){(\emph{b})}
\put(55,15){\circle*{1.0}}
\put(55,30){\circle{1.0}}
\put(70,15){\circle{1.0}}
\put(70,30){\circle{1.0}}
\put(55,15){\vector(0,1){15.00}}
\put(55,15){\vector(1,0){15.00}}
\put(55,15){\line(1,1){15.00}}

\bezier{20}(55.5,15)(55.5,22.5)(55.5,30)
\bezier{20}(56.0,15)(56.0,22.5)(56.0,30)
\bezier{20}(56.5,15)(56.5,22.5)(56.5,30)
\bezier{20}(57.0,15)(57.0,22.5)(57.0,30)
\bezier{20}(57.5,15)(57.5,22.5)(57.5,30)
\bezier{20}(58.0,15)(58.0,22.5)(58.0,30)
\bezier{20}(58.5,15)(58.5,22.5)(58.5,30)
\bezier{20}(59.0,15)(59.0,22.5)(59.0,30)
\bezier{20}(59.5,15)(59.5,22.5)(59.5,30)
\bezier{20}(60.0,15)(60.0,22.5)(60.0,30)
\bezier{20}(60.5,15)(60.5,22.5)(60.5,30)
\bezier{20}(61.0,15)(61.0,22.5)(61.0,30)
\bezier{20}(61.5,15)(61.5,22.5)(61.5,30)
\bezier{20}(62.0,15)(62.0,22.5)(62.0,30)
\bezier{20}(62.5,15)(62.5,22.5)(62.5,30)
\bezier{20}(63.0,15)(63.0,22.5)(63.0,30)
\bezier{20}(63.5,15)(63.5,22.5)(63.5,30)
\bezier{20}(64.0,15)(64.0,22.5)(64.0,30)
\bezier{20}(64.5,15)(64.5,22.5)(64.5,30)
\bezier{20}(65.0,15)(65.0,22.5)(65.0,30)
\bezier{20}(65.5,15)(65.5,22.5)(65.5,30)
\bezier{20}(66.0,15)(66.0,22.5)(66.0,30)
\bezier{20}(66.5,15)(66.5,22.5)(66.5,30)
\bezier{20}(67.0,15)(67.0,22.5)(67.0,30)
\bezier{20}(67.5,15)(67.5,22.5)(67.5,30)
\bezier{20}(68.0,15)(68.0,22.5)(68.0,30)
\bezier{20}(68.5,15)(68.5,22.5)(68.5,30)
\bezier{20}(69.0,15)(69.0,22.5)(69.0,30)
\bezier{20}(69.5,15)(69.5,22.5)(69.5,30)
\bezier{20}(70,15)(70,22.5)(70,30)
\put(53,12){$v$}
\put(70,12){$v$}
\put(53,31){$v$}
\put(70,31){$v$}
\put(61,12){$a_1$}
\put(51,24){$a_2$}
\put(65,17){$\Omega$}

\put(98.0,15){\circle*{1}}
\put(97,12){$v$}
\put(103,20){$\scriptstyle(1,0)$}
\put(85,20){$\scriptstyle(0,1)$}
\put(94,31){$\scriptstyle(1,1)$}
\bezier{200}(98.0,15)(85,22)(84.0,15)
\bezier{200}(98.0,15)(85,8)(84.0,15)
\bezier{200}(98.0,15)(111,22)(112.0,15)
\bezier{200}(98.0,15)(111,8)(112.0,15)
\bezier{200}(98.0,15)(105,28)(98.0,29)
\bezier{200}(98.0,15)(91,28)(98.0,29)
\end{picture}
\caption{\footnotesize \emph{a}) The triangular lattice, the fundamental cell $\Omega$ is shaded; \emph{b})~its fundamental graph with the minimal form $\gm\in\mF(\k)$.}
\lb{ff.10T}
\end{figure}
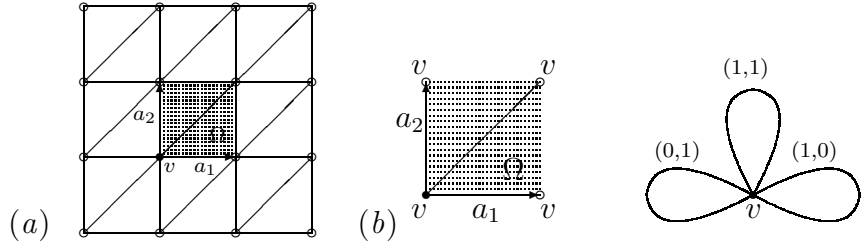

\medskip

\no \textbf{Proof of Corollary \ref{TCo3}.} Theorem \ref{Pphi}.\emph{i}) gives that the minimal form $\gm$ satisfies $\gm=\wt\gm(\,\cdot\,,T)$ for some $\gm$-minimal spanning tree $T$ of the graph $G$, where $\wt\gm(\,\cdot\,,T):A\ra\R^d$ is defined by \er{cat} as $\mathbf{x}=\gm$. Since $\Phi_{\gm}(\cC)=\Z^d$, we obtain $\gm(A)\ss\Z^d$. We construct a $\Z^d$-periodic graph $\cG$ such that $G$ is its fundamental graph and $\gm$ is the index form on $G$. We identify all vertices of the graph $G$ with some distinct points of $[0,1)^d$. Each edge $\be=(u,v)$ of the graph $G$ induces an infinite number of edges
$$
\big\{\big(u+\mn,v+\mn+\gm(\be)\big)\big\}_{\mn\in\Z^d}
$$
of the graph $\cG$ with the edge index $\gm(\be)$. By construction, $\cG$ is a $\Z^d$-periodic graph with the fundamental graph $G_*=G$ and $\gm$ is the index form on $G_*$. Thus, due to Theorem \ref{TFD1}, the operator $\mA(\vt)$, $\vt\in\T^d$, defined by \er{l2.13at} is a fiber operator for the Laplacian on the graph $\cG$. \qq \BBox

\section{\lb{Sec4} Proof of the estimates}
\setcounter{equation}{0}
\subsection{The Lebesgue measure of the spectrum}
In order to prove Theorem \ref{T1} we need the following lemma.

\begin{lemma}\lb{fEst} Let $\vt\in\T^d$, $f\in\ell^2(\cV_*)$, $\gm\in\mF(\k)$ be a minimal form on the fundamental graph $G_*=(\cV_*,\cE_*)$, where $\mF(\k)$ is given by \er{vv1f} at $\mathbf{x}=\k$, and $\k$ is the coordinate form defined by \er{edco}, \er{dco}. Then

i) The fiber Laplacian $\D_\gm(\vt)$ given by \er{l2.13am} has the following representation:
\[\lb{fDmDb}
\D_\gm(\vt)=\D_0+\wt\D_\gm(\vt),
\]
where $\D_0$ is the Laplacian on the finite graph $G_\gm^0=(\cV_*,\cE_*\sm\cE_\gm)$, $\cE_\gm=\cE_*\cap\supp \gm$,
$\wt\D_\gm(\vt)$ is the magnetic Laplacian with the magnetic vector potential $\a(\be)=\lan\gm(\be),\vt\ran$, $\be\in\cA_\gm$, on the finite graph
$G_\gm=(\cV_*,\cE_\gm)$:
\[\label{fDbt}
\big(\wt\D_\gm(\vt)f\big)(v)=\sum_{\be=(v,\,u)\in\cA_\gm}
\big(f(v)-e^{i\lan\gm(\be),\,\vt\ran}f(u)\big), \qq
v\in\cV_*,
\]
$\cA_\gm=\{\be \mid \be\in\cE_\gm \; \textrm{or} \;\; \ul\be\in\cE_\gm\}$.

ii) The quadratic form of the magnetic Laplacian
$\wt\D_\gm(\vt)$ is given by
\[\lb{qflom}
\lan \wt\D_\gm(\vt)
f,f\ran_{\ell^2(\cV_*)}={1\/2}\sum_{\be=(v,u)\in\cA_\gm}
\big|f(v)-e^{i\,\lan\gm(\be),\,\vt\ran}f(u)\big|^2.
\]

iii) The magnetic Laplacian $\wt\D_\gm(\vt)$ given by \er{fDbt} satisfies
\[\lb{0DaB}
0\le\wt\D_\gm(\vt)\le2B_\gm,
\]
where
\[\lb{Bfv}
(B_\gm f)(v)=\vk_v^{\gm}f(v), \qqq v\in\cV_*,
\]
$\vk_v^{\gm}$ is the degree of the vertex $v\in\cV_*$ on the graph $G_\gm=(\cV_*,\cE_\gm)$, and
\[\lb{fsdbe}
\s\big(\wt\D_\gm(\vt)\big)\ss[0,2\vk_+^{\gm}], \qqq \vk^{\gm}_+=\max_{v\in\cV_*}\vk_v^{\gm}.
\]
\end{lemma}

\no \textbf{Proof.} \emph{i}) Let $f\in\ell^2(\cV_*)$, $v\in\cV_*$. Then using \er{l2.13am} for each $\vt\in\T^d$ we obtain
$$
\begin{aligned}
&\big(\D_\gm(\vt)f\big)(v)=\sum_{\be=(v,\,u)\in\cA_*}
\big(f(v)-e^{i\lan\gm(\be),\,\vt\ran}f(u)\big)=
\sum_{\be=(v,\,u)\in\cA_*\sm\cA_\gm}\big(f(v)-f(u)\big)\\&+
\sum_{\be=(v,\,u)\in\cA_\gm}\big(f(v)-e^{i\lan\gm(\be),\,\vt\ran}f(u)\big)=
(\D_0f)(v)+\big(\wt\D_\gm(\vt)f\big)(v).
\end{aligned}
$$
For each $\vt\in\T^d$ the operator $\wt\D_\gm(\vt)$ is the magnetic Laplacian with the magnetic potential $\a(\be)=\lan\gm(\be),\vt\ran$, $\be\in\cA_\gm$, on the finite graph $G_\gm=(\cV_*,\cE_\gm)$.

\emph{ii}) This item was proved in \cite{KS17} (see Theorem 6.4.iv).

\emph{iii}) We show that $2B_\gm-\wt\D_\gm(\vt)\ge0$ for each $\vt\in\T^d$. Let $f\in\ell^2(\cV_*)$. Then, using the identities \er{qflom}, \er{Bfv}, we obtain
\[\lb{2BD0}
\begin{aligned}
&\textstyle\lan\big(2B_\gm-\wt\D_\gm(\vt)\big)f,f\ran_{\ell^2(\cV_*)}=
2\sum\limits_{v\in\cV_*}\vk_v^\gm|f(v)|^2-\lan\wt\D_\gm(\vt)f,f\ran_{\ell^2(\cV_*)}\\
&\textstyle=\sum\limits_{(v,\,u)\in\cA_\gm}\big(|f(v)|^2+
|f(u)|^2\big)-{1\/2}\sum\limits_{\be=(v,\,u)\in\cA_\gm}
\big|f(v)-e^{i\lan\gm(\be),\,\vt\ran}f(u)\big|^2\\
&\textstyle={1\/2}\sum\limits_{\be=(v,\,u)\in\cA_\gm}
\big|f(v)+e^{i\lan\gm(\be),\,\vt\ran}f(u)\big|^2\geq0.
\end{aligned}
\]
For each\, $\vt\in\T^d$ the spectrum of the magnetic operator $\wt\D_\gm(\vt)$  satisfies the condition \er{fsdbe} (see, e.g., \cite{HS99a}), which, in particular, yields that $\wt\D_\gm(\vt)\geq0$. \qq \BBox

\medskip

\no {\bf Proof of Theorem \ref{T1}.} \emph{i}) Using \er{fDmDb}, we rewrite the fiber Schr\"odinger operator $H_\gm(\vt)$,  $\vt\in\T^d$, defined by \er{Hvtm}, \er{l2.13am} in the form:
\[
\label{eq.1}
H_\gm(\vt)=H_0+\wt\D_\gm(\vt),
\]
where $H_0=\D_0+Q$ is the Schr\"odinger operator on the finite graph $G_\gm^0=(\cV_*,\cE_*\sm\cE_\gm)$, $\wt\D_\gm(\vt)$ is the magnetic Laplacian on the finite graph $G_\gm=(\cV_*,\cE_\gm)$ defined by \er{fDbt}. Lemma \ref{fEst}.\emph{iii} gives that the spectrum $\s\big(\wt\D_\gm(\vt)\big)\ss[0,2\vk^\gm_+]$, where $\vk^\gm_+$ is defined in \er{fesbp1}. Then each eigenvalue $\l_n(\vt)$, $n\in\N_\n$, of $H_\gm(\vt)$ satisfies $\m_{\gm,n}\le\l_n(\vt)\le\m_{\gm,n}+2\vk^\gm_+$, which yields
\er{fesbp1}.

\emph{ii}) From \er{eq.1} and \er{0DaB} we obtain
$$
H_0\le H_\gm(\vt)\le H_0+2B_\gm,
$$
where $B_\gm$ is defined in \er{Bfv}. This yields
$$
\l_n(H_0)\leq\l_n^-\le \l_n(\vt)\leq\l_n^+\le \l_n(H_0+2B_\gm),
\qqq \forall\,(n,\vt)\in\N_\n\ts\T^d.
$$
Then, using \er{dIm},
\begin{multline*}
\big|\s(H)\big|\le\sum_{n=1}^{\nu}(\l_n^+-\l_n^-)\leq
\sum_{n=1}^\nu\big(\l_n(H_0+2B_\gm)-\l_n(H_0)\big)\\=2\Tr B_\gm=2\sum_{v\in\cV_*}\vk^\gm_v=
4\#\cE_\gm=2\#\supp \gm=4\cI.
\end{multline*}

\emph{iii}) This item follows from Proposition \ref{TG1}.
\qq $\BBox$

\setlength{\unitlength}{1.0mm}
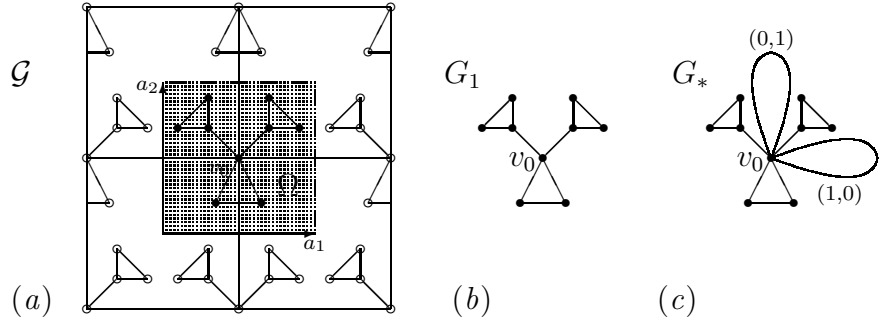
\begin{figure}[h]
\centering
\unitlength 1mm 
\linethickness{0.4pt}
\ifx\plotpoint\undefined\newsavebox{\plotpoint}\fi 
\begin{picture}(120,50)(0,0)
\put(0,40){$\cG$}
\put(0,10){(\emph{a})}
\put(20,20){\vector(1,0){20.00}}
\put(20,20){\vector(0,1){20.00}}
\multiput(21,40)(4,0){5}{\line(1,0){2}}
\multiput(40,21)(0,4){5}{\line(0,1){2}}
\put(38.5,18.0){$\scriptstyle a_1$}
\put(16.5,39.0){$\scriptstyle a_2$}
\put(26,28){$\scriptstyle v_0$}
\put(35,25){$\Omega$}

\bezier{40}(20.5,20)(20.5,30)(20.5,40)
\bezier{40}(21.0,20)(21.0,30)(21.0,40)
\bezier{40}(21.5,20)(21.5,30)(21.5,40)
\bezier{40}(22.0,20)(22.0,30)(22.0,40)
\bezier{40}(22.5,20)(22.5,30)(22.5,40)
\bezier{40}(23.0,20)(23.0,30)(23.0,40)
\bezier{40}(23.5,20)(23.5,30)(23.5,40)
\bezier{40}(24.0,20)(24.0,30)(24.0,40)
\bezier{40}(24.5,20)(24.5,30)(24.5,40)
\bezier{40}(25.0,20)(25.0,30)(25.0,40)
\bezier{40}(25.5,20)(25.5,30)(25.5,40)
\bezier{40}(26.0,20)(26.0,30)(26.0,40)
\bezier{40}(26.5,20)(26.5,30)(26.5,40)
\bezier{40}(27.0,20)(27.0,30)(27.0,40)
\bezier{40}(27.5,20)(27.5,30)(27.5,40)
\bezier{40}(28.0,20)(28.0,30)(28.0,40)
\bezier{40}(28.5,20)(28.5,30)(28.5,40)
\bezier{40}(29.0,20)(29.0,30)(29.0,40)
\bezier{40}(29.5,20)(29.5,30)(29.5,40)
\bezier{40}(30.0,20)(30.0,30)(30.0,40)

\bezier{40}(30.5,20)(30.5,30)(30.5,40)
\bezier{40}(31.0,20)(31.0,30)(31.0,40)
\bezier{40}(31.5,20)(31.5,30)(31.5,40)
\bezier{40}(32.0,20)(32.0,30)(32.0,40)
\bezier{40}(32.5,20)(32.5,30)(32.5,40)
\bezier{40}(33.0,20)(33.0,30)(33.0,40)
\bezier{40}(33.5,20)(33.5,30)(33.5,40)
\bezier{40}(34.0,20)(34.0,30)(34.0,40)
\bezier{40}(34.5,20)(34.5,30)(34.5,40)
\bezier{40}(35.0,20)(35.0,30)(35.0,40)
\bezier{40}(35.5,20)(35.5,30)(35.5,40)
\bezier{40}(36.0,20)(36.0,30)(36.0,40)
\bezier{40}(36.5,20)(36.5,30)(36.5,40)
\bezier{40}(37.0,20)(37.0,30)(37.0,40)
\bezier{40}(37.5,20)(37.5,30)(37.5,40)
\bezier{40}(38.0,20)(38.0,30)(38.0,40)
\bezier{40}(38.5,20)(38.5,30)(38.5,40)
\bezier{40}(39.0,20)(39.0,30)(39.0,40)
\bezier{40}(39.5,20)(39.5,30)(39.5,40)
\bezier{40}(40.0,20)(40.0,30)(40.0,40)

\put(10,10){\line(1,0){40.00}}
\put(10,30){\line(1,0){40.00}}
\put(30,10){\line(0,1){40.00}}
\put(10,50){\line(1,0){40.00}}
\put(10,10){\line(0,1){40.00}}
\put(50,10){\line(0,1){40.00}}

\put(10,10){\circle{1}}
\put(30,10){\circle{1}}
\put(50,10){\circle{1}}

\put(10,30){\circle{1}}
\put(30,30){\circle*{1}}
\put(50,30){\circle{1}}

\put(10,50){\circle{1}}
\put(30,50){\circle{1}}
\put(50,50){\circle{1}}
\put(30,30){\line(1,1){4.00}}
\put(34,34){\line(1,0){4.00}}
\put(34,34){\line(0,1){4.00}}
\put(34,38){\line(1,-1){4.00}}
\put(34,34){\circle*{1}}
\put(34,38){\circle*{1}}
\put(38,34){\circle*{1}}

\put(30,30){\line(-1,1){4.00}}
\put(26,34){\line(-1,0){4.00}}
\put(26,34){\line(0,1){4.00}}
\put(26,38){\line(-1,-1){4.00}}
\put(26,34){\circle*{1}}
\put(26,38){\circle*{1}}
\put(22,34){\circle*{1}}

\put(30,30){\line(-1,-2){3.00}}
\put(30,30){\line(1,-2){3.00}}
\put(27,24){\line(1,0){6.00}}
\put(27,24){\circle*{1}}
\put(33,24){\circle*{1}}
\put(65.5,29){$v_0$}
\put(50,30){\line(-1,1){4.00}}
\put(46,34){\line(-1,0){4.00}}
\put(46,34){\line(0,1){4.00}}
\put(46,38){\line(-1,-1){4.00}}
\put(46,34){\circle{1}}
\put(46,38){\circle{1}}
\put(42,34){\circle{1}}

\put(50,30){\line(-1,-2){3.00}}
\put(47,24){\line(1,0){3.00}}
\put(47,24){\circle{1}}
\put(10,30){\line(1,1){4.00}}
\put(14,34){\line(1,0){4.00}}
\put(14,34){\line(0,1){4.00}}
\put(14,38){\line(1,-1){4.00}}
\put(14,34){\circle{1}}
\put(14,38){\circle{1}}
\put(18,34){\circle{1}}

\put(10,30){\line(1,-2){3.00}}
\put(10,24){\line(1,0){3.00}}
\put(13,24){\circle{1}}
\put(30,10){\line(1,1){4.00}}
\put(34,14){\line(1,0){4.00}}
\put(34,14){\line(0,1){4.00}}
\put(34,18){\line(1,-1){4.00}}
\put(34,14){\circle{1}}
\put(34,18){\circle{1}}
\put(38,14){\circle{1}}

\put(30,10){\line(-1,1){4.00}}
\put(26,14){\line(-1,0){4.00}}
\put(26,14){\line(0,1){4.00}}
\put(26,18){\line(-1,-1){4.00}}
\put(26,14){\circle{1}}
\put(26,18){\circle{1}}
\put(22,14){\circle{1}}

\put(50,10){\line(-1,1){4.00}}
\put(46,14){\line(-1,0){4.00}}
\put(46,14){\line(0,1){4.00}}
\put(46,18){\line(-1,-1){4.00}}
\put(46,14){\circle{1}}
\put(46,18){\circle{1}}
\put(42,14){\circle{1}}

\put(10,10){\line(1,1){4.00}}
\put(14,14){\line(1,0){4.00}}
\put(14,14){\line(0,1){4.00}}
\put(14,18){\line(1,-1){4.00}}
\put(14,14){\circle{1}}
\put(14,18){\circle{1}}
\put(18,14){\circle{1}}
\put(30,50){\line(-1,-2){3.00}}
\put(30,50){\line(1,-2){3.00}}
\put(27,44){\line(1,0){6.00}}
\put(27,44){\circle{1}}
\put(33,44){\circle{1}}
\put(50,50){\line(-1,-2){3.00}}
\put(47,44){\line(1,0){3.00}}
\put(47,44){\circle{1}}
\put(10,50){\line(1,-2){3.00}}
\put(10,44){\line(1,0){3.00}}
\put(13,44){\circle{1}}

\put(58,10){(\emph{b})}
\put(100,30){\circle*{1}}

\put(100,30){\line(1,1){4.00}}
\put(104,34){\line(1,0){4.00}}
\put(104,34){\line(0,1){4.00}}

\put(104,38){\line(1,-1){4.00}}
\put(104,34){\circle*{1}}
\put(104,38){\circle*{1}}
\put(108,34){\circle*{1}}

\put(100,30){\line(-1,1){4.00}}
\put(96,34){\line(-1,0){4.00}}
\put(96,34){\line(0,1){4.00}}

\put(96,38){\line(-1,-1){4.00}}
\put(96,34){\circle*{1}}
\put(96,38){\circle*{1}}
\put(92,34){\circle*{1}}

\put(100,30){\line(-1,-2){3.00}}
\put(100,30){\line(1,-2){3.00}}

\put(97,24){\circle*{1}}
\put(103,24){\circle*{1}}
\put(97,24){\line(1,0){6.00}}

\bezier{300}(100,30)(105,43)(100,44)
\bezier{300}(100,30)(95,43)(100,44)

\bezier{300}(100,30)(113,35)(114,30)
\bezier{300}(100,30)(113,25)(114,30)
\put(95.5,29){$v_0$}
\put(87,40){$G_*$}
\put(97,45){$\scriptstyle(0,1)$}
\put(106,24.5){$\scriptstyle(1,0)$}
\put(85,10){(\emph{c}) }

\put(57,40){$G_1$}
\put(70,30){\circle*{1}}
\put(70,30){\line(1,1){4.00}}
\put(74,34){\line(1,0){4.00}}
\put(74,34){\line(0,1){4.00}}
\put(74,38){\line(1,-1){4.00}}
\put(74,34){\circle*{1}}
\put(74,38){\circle*{1}}
\put(78,34){\circle*{1}}

\put(70,30){\line(-1,1){4.00}}
\put(66,34){\line(-1,0){4.00}}
\put(66,34){\line(0,1){4.00}}
\put(66,38){\line(-1,-1){4.00}}
\put(66,34){\circle*{1}}
\put(66,38){\circle*{1}}
\put(62,34){\circle*{1}}

\put(70,30){\line(-1,-2){3.00}}
\put(70,30){\line(1,-2){3.00}}
\put(67,24){\line(1,0){6.00}}
\put(67,24){\circle*{1}}
\put(73,24){\circle*{1}}

\end{picture}
\vspace{-0.5cm}
\caption{\footnotesize  \emph{a}) a $\G$-periodic graph $\cG$, $\{a_1,a_2\}$ is the basis of the lattice $\G$; \; \emph{b}) a finite decoration $G_1$;\quad \emph{c}) the fundamental graph $G_*$, the values of the minimal form $\gm(\be)$, $\be\in\supp\gm$, are shown near the edges.}
\label{ff.10}
\end{figure}

\begin{proposition}\lb{TG1}
Let $\cG$ be a periodic graph obtained from the
$d$-dimensional lattice by "gluing"\, the same finite connected graph $G_1$ with the Betti number $\b_1$ to each vertex of the lattice (for $d=2$ see Fig.\ref{ff.10}a). Then the following statements hold true.

i) For the periodic graph $\cG$, the Betti number $\b$, the invariant $\cI$ given by \er{dIm} and the maximum vertex degree $\vk_+$ defined in \er{bf} satisfy
\[\lb{deco}
\b=\b_1+d, \qqq \cI=d, \qqq \vk_+-2\cI\geq\vk_{v_0}^{(1)},
\]
where $\vk_{v_0}^{(1)}$ is the degree of the "gluing" vertex $v_0$ on the graph $G_1$.

ii) The spectrum of the Schr\"odinger operator $H=\D+Q$ on $\cG$ has the form
\[
\lb{acs3}
\s(H)=\bigcup_{n=1}^\n\big[\l_n(0),\l_n(\ol\pi)\big], \qqq \ol\pi=(\pi,\ldots,\pi)\in\T^d.
\]

iii) The Lebesgue measure of the spectrum of the Schr\"odinger operators $H$ on $\cG$ satisfies
\[
\lb{sp2}
\textstyle|\s(H)|=4d=4\cI
\]
and the estimate \er{eq.7'} becomes an identity.
\end{proposition}

\no {\bf Proof.} \emph{i}) The fundamental graph $G_*=(\cV_*,\cE_*)$ of the periodic graph $\cG$ consists of all edges of the graph $G_1$ and $d$ loops $\be_1,\ldots,\be_d$ of the fundamental graph of the $d$-dimensional lattice. Then $\vk_{v_0}=2d+\vk_{v_0}^{(1)}$ and $\b=\b_1+d$. Using Theorem \ref{Pphi}.\emph{i}) -- \emph{ii}) we obtain that the minimal form $\gm\in\mF(\k)$ is unique and $\supp\gm=\{\be_1,\ldots,\be_d,\ul\be_1,\ldots,\ul\be_d\}$, which yields the identity $\cI=d$. Then, using the definition \er{bf} of $\vk_+$, we obtain
$$
\vk_+=\max_{v\in\cV}\vk_v\geq\vk_{v_0}=2d+\vk_{v_0}^{(1)}=2\cI+\vk_{v_0}^{(1)}.
$$

\emph{ii}) -- \emph{iii}) The identity \er{acs3} and the first identity in \er{sp2} were proved in \cite{KS14} (Proposition 7.2). Since $d=\cI$, for the considered graph $\cG$ the estimate \er{eq.7'} becomes an identity. \qq
\BBox

\subsection{Effective masses at the bottom of the spectrum}
We introduce the Hilbert space
\[
\ell^2(\cA_*)=\{\phi: \cA_*\ra\C\mid \phi(\ul\be)=-\phi(\be)
\textrm{ for } \be\in\cA_* \qq \textrm{ and } \qq \lan
\phi,\phi\ran_{\cA_*}<\iy\},
\]
where the inner product is given by
\[\lb{incA}
\lan
\phi_1,\phi_2\ran_{\cA_*}={1\/2}\sum_{\be\in\cA_*}\phi_1(\be)\ol{\phi_2(\be)}.
\]

Let $\gm\in\mF(\k)$ be a minimal form on the fundamental graph $G_*=(\cV_*,\cE_*)$, where $\mF(\k)$ is given by \er{vv1f} at $\mathbf{x}=\k$, and $\k$ is the coordinate form defined by \er{edco}, \er{dco}. For each $\vt\in\T^d$ we define the operator
$\na_\gm(\vt):\ell^2(\cV_*)\ra \ell^2(\cA_*)$ by
\[\lb{navt}
\big(\na_\gm(\vt)
f\big)(\be)=e^{-i\lan\gm(\be),\vt\ran/2}f(v)-e^{i\lan\gm(\be),\vt\ran/2}f(u),\qq \be=(v,u), \qqq \forall\,f\in \ell^2(\cV_*).
\]
We need the following results about the factorization of fiber Laplacians (see Theorem 6.4.i -- ii) in \cite{KS17}).

\emph{i) For each $\vt\in\T^d$ the conjugate operator
$\na^*_\gm(\vt):\ell^2(\cA_*)\ra \ell^2(\cV_*)$ has
the form}
\[\lb{coop}
(\na^*_\gm(\vt)
\phi)(v)=\sum_{\be=(v,u)\in\cA_*}e^{i\lan\gm(\be),\vt\ran/2}\phi(\be),\qqq
\forall\,\phi\in \ell^2(\cA_*).
\]

\emph{ii) For each $\vt\in\T^d$ the fiber Laplacian
$\D_\gm(\vt)$ defined by \er{l2.13am} satisfies
\[\lb{fact}
\D_\gm(\vt)=\na^*_\gm(\vt)\na_\gm(\vt).
\]
}

The Taylor expansion of $\na_\gm(\vt)$ about the point $\vt_0=0$ is given by
\[
\lb{nas}
\na_\gm(\vt)=\na_\gm(0)+\ve\na_1(\o)+\ve^2\na_2(\o)+O(\ve^3) \qq
\textrm{as } \ \vt=\ve\o, \qqq \ve=|\vt\,|\rightarrow0,
\]
where
\[
\lb{anb}
\na_1(\o)={\dot\na_\gm(\ve\o)}\big|_{\ve=0},\qqq
\na_2(\o)=\textstyle{1\/2}\,\ddot\na_\gm(\ve\o)\big|_{\ve=0}, \qqq\o\in \S^{d-1}\,,
\]
$\dot u=\pa u/\pa \ve$ and $\S^{d}$ is the $d$-dimensional sphere.

\medskip

Let $\l(\vt)$, $\vt\in\T^d$, be the first band function of the Laplacian $\D$. It is known that $\l(0)=0$ is a simple eigenvalue of $\D_\gm(0)$ with the eigenfunction
\[\lb{psi0}
\p(0,\cdot)=\1_\n=(1,\ldots,1)\in\R^\n.
\]
Then the eigenvalue $\l(\vt)$ of $\D_\gm(\vt)$ and the
corresponding eigenfunction $\p(\vt,\cdot)$ have asymptotics as $\vt=\ve\o$, $\o\in\S^{d-1}$, $\ve\to0$:
\[
\label{lamp}
\begin{aligned}
\l(\vt)=\ve^2\m(\o)+O(\ve^3), \qqq
\p(\vt,\cdot)=\1_\n+\ve\p_1+\ve^2\p_2+O(\ve^3),\\
\textstyle\m(\o)={1\/2}\,\ddot\l(\ve\o)\big|_{\ve=0},\qq
\p_1=\p_1(\o,\cdot)=\dot\p(\ve\o,\cdot)\big|_{\ve=0},\qq
\textstyle\p_2=\p_2(\o,\cdot)={1\/2}\,\ddot\p(\ve\o,\cdot)\big|_{\ve=0}.
\end{aligned}
\]

\begin{theorem}
\label{Tem}  Let $\gm\in\mF(\k)$ be a minimal form on the fundamental graph $G_*=(\cV_*,\cE_*)$. Then the effective form $\mu(\o)$, $\o\in\S^{d-1}$, defined in \er{lamp} satisfies
\[
\lb{tseN}
{1\/\n^2d}\leq\mu(\o)\leq{1\/2\n}\sum\limits_{\be\in\supp\gm}\lan\gm
(\be),\o\ran^2,\qqq  \n=\#\cV_*.
\]
\end{theorem}

\no{\bf Proof.} The fiber Laplacian $\D_\gm(\vt)$, $\vt\in\T^d$, defined by \er{l2.13am}, can be represented in the following form:
\[
\label{geq.1}
\D_\gm(\vt)=\D_\gm(0)+\ve\D_1(\o)+\ve^2\D_2(\o)+O(\ve^3), \\
\]
as $\vt=\ve\o$, $\ve\rightarrow0$, $\o\in\S^{d-1}$, where
\[\lb{DDD}
\D_1(\o)=\dot
\D_\gm(\ve\o)\big|_{\ve=0},\qqq \D_2(\o)=\textstyle{1\/2}\,\ddot
\D_\gm(\ve\o)\big|_{\ve=0}.
\]
The equation $\D_\gm(\vt)\p(\vt,\cdot)=\l(\vt)\p(\vt,\cdot)$ after substitution
\er{lamp}, \er{geq.1} takes the form
\[\lb{gaas}
\begin{aligned}
\big(\D_\gm(0)+\ve\D_1(\o)+\ve^2\D_2(\o)+O(\ve^3)\big)
\big(\1_\n+\ve\p_1+\ve^2\p_2+O(\ve^3)\big)\\=
\big(\ve^2\m(\o)+O(\ve^3)\big)
\big(\1_\n+\ve\p_1+\ve^2\p_2+O(\ve^3)\big),
\end{aligned}
\]
where $\p_1,\p_2$ are defined in \er{lamp}. This asymptotics gives two identities for any $\o\in \S^{d-1}$:
\[\lb{gve1}
\D_1(\o)\1_\n+\D_\gm(0)\p_1=0,
\]
\[\lb{gve2}
\D_2(\o)\1_\n+\D_1(\o)\p_1+\D_\gm(0)\p_2=\m(\o)\,\1_\n.
\]

Substituting \er{nas} and \er{geq.1} into the identity \er{fact} we obtain
\[\lb{Dna}
\begin{aligned}
&\D_1(\o)=\na_\gm^*(0)\na_1(\o)+\na_1^*(\o)\na_\gm(0),\\
&\D_2(\o)=\na_\gm^*(0)\na_2(\o)+\na_1^\ast(\o)\na_1(\o)+\na_2^*(\o)\na_\gm(0),
\end{aligned}
\]
where $\D_s(\o)$ and $\na_s(\o)$, $s=1,2$, are defined in \er{DDD} and \er{anb}, respectively.

Multiplying both sides of \er{gve1} by $\p_1$ and both sides of
\er{gve2} by $\1_\n$ and using the identities \er{fact}, \er{Dna} and $\na_\gm(0)\1_\n=0$, we have
\[\lb{ve1'}
\begin{aligned}
&\lan\na_1(\o)\1_\n,\na_\gm(0)\p_1\ran_{\cA_*}+\|\na_\gm(0)\p_1\|^2_{\cA_*}=0,\\
&\|\na_1(\o)\1_\n\|^2_{\cA_*}+\lan\na_\gm(0)\p_1,\na_1(\o)\1_\n\ran_{\cA_*}=\n\m(\o).
\end{aligned}
\]
Then, combining \er{ve1'}, we obtain
\[\lb{mu2}
\n\m(\o)=\|\na_1(\o)\1_\n\|^2_{\cA_*}-\|\na_\gm(0)\p_1\|_{\cA_*}^2=
\|\na_1(\o)\1_\n+\na_\gm(0)\p_1\|_{\cA_*}^2.
\]
In virtue of Theorem \ref{Pphi}.\emph{i}), the minimal form $\gm$ satisfies $\gm=\wt\gm(\,\cdot\,,T)$ for some $\k$-minimal spanning tree $T=(\cV_*,\cE_T)$ of the fundamental graph $G_*$, where $\wt\gm(\,\cdot\,,T): \cA_*\ra\R^d$ is defined by \er{cat} as $\mathbf{x}=\k$. Let $\be\in\supp\gm$. Then, due to the property 2) of spanning trees (see page \pageref{PSTs}), there exists a unique cycle $\mathbf{c}_{\be}$ whose edges are all in $T$ except $\be$. Let $\phi_\be\in\ell^2(\cA_*)$ be defined by
\[\lb{xsi}
\phi_\be(\be_*)=\left\{
\begin{array}{rc}
  1, & \textrm{if }  \be_*\in\mathbf{c}_{\be}\\
  -1, & \textrm{if }  \ul\be_*\in\mathbf{c}_{\be}\\
  0, & \textrm{otherwise}
\end{array}\right..
\]
Since the length of $\mathbf{c}_{\be}$ is not more than $\n$, $\|\phi_\be\|_{\cA_*}^2\leq\n$. Then, using \er{mu2} and the Cauchy-Schwarz inequality, we have
\begin{multline}\lb{mmo}
\n^2\m(\o)\geq\n\m(\o)\,\|\phi_\be\|_{\cA_*}^2=
\|\na_1(\o)\1_\n+\na_\gm(0)\p_1\|_{\cA_*}^2\cdot\|\phi_\be\|_{\cA_*}^2\geq
|\lan\na_1(\o)\1_\n+\na_\gm(0)\p_1,\phi_\be\ran_{\cA_*}|^2\\=
|\lan\na_1(\o)\1_\n,\phi_\be\ran_{\cA_*}+\lan\na_\gm(0)\p_1,\phi_\be\ran_{\cA_*}|^2=
|\lan\na_1(\o)\1_\n,\phi_\be\ran_{\cA_*}+\lan\p_1,\na_\gm^*(0)\phi_\be\ran_{\cA_*}|^2.
\end{multline}
Using the first identity in \er{anb} and \er{navt} we have
\[
\label{nab1}
\big(\na_1(\o)\1_\n\big)(\be)=-i\,\lan\gm(\be),\o\ran.
\]
Due to \er{coop}, $\na_\gm^*(0)\phi_\be=0$ for each $\be\in\supp\gm$. Substituting this and \er{nab1} into \er{mmo} and using the identity
$\Phi_{\gm}(\mathbf{c}_{\be})=\gm(\be)$, we obtain
\begin{multline*}
\n^2\m(\o)\geq\big|\lan\na_1(\o)\1_\n,\phi_\be\ran_{\cA_*}\big|^2=
\big|\textstyle\frac12\sum\limits_{\wt\be\in\cA_*}\phi_\be(\wt\be\,)
\lan\gm(\wt\be\,),\o\ran\big|^2\\
=\big|\textstyle\frac12
\big\lan\sum\limits_{\wt\be\in\cA_*}\phi_\be(\wt\be\,)\gm(\wt\be\,),\o\big\ran\big|^2=
\big\lan\sum\limits_{\wt\be\in\mathbf{c}_{\be}}\gm(\wt\be\,),\o\big\ran^2=\lan
\gm(\be),\o\ran^2.
\end{multline*}
Thus, we have
\[\lb{esmu}
\n^2\m(\o)\geq\lan
\gm(\be),\o\ran^2,\qqq \forall\, \be\in\supp\gm.
\]
Due to Theorem \ref{TDImf}.\emph{i}), for some basis $a_1,\ldots,a_d$ of the lattice $\G$ there exist edges $\be_1,\ldots,\be_d\in\supp\gm$ such that the set
$\big\{\gm(\be_s)\big\}_{s\in\N_d}$ forms an orthonormal basis of $\Z^d$. Summing the estimates \er{esmu} over all these edges $\be_1,\ldots,\be_d$, we get
\[\lb{esm1}
d\,\n^2\m(\o)\geq\textstyle\sum\limits_{s=1}^d\lan
\gm(\be_s),\o\ran^2=|\o|^2=1.
\]
Thus, the lower bound in \er{tseN} has been proved.

Now we prove the upper bound in \er{tseN}. From \er{mu2} we deduce
\[\lb{upe}
\n\m(\o)\leq\|\na_1(\o)\1_\n\|^2_{\cA_*}.
\]
Substituting \er{nab1} into \er{upe}, we obtain
\[
\lb{eem}
\n\m(\o)\leq\textstyle\frac12\sum\limits_{\be\in\cA_*}\lan\gm(\be),\o\ran^2=
\textstyle\frac12\sum\limits_{\be\in\supp\gm}\lan\gm(\be),\o\ran^2.
\]
\BBox

\medskip

\no \textbf{Proof of Theorem \ref{TTsem}.} Let $0<M_1\le M_2\le\ldots\le M_d$ and $\wt\o_1,\ldots,\wt\o_d$ be the eigenvalues and the corresponding orthonormal eigenvectors of the matrix $M$. From the definition of the effective form $\mu(\o)$, $\o\in\S^{d-1}$, in \er{lamp} and \er{lam} it follows that $\mu(\o)=\frac12\lan M\o,\o\ran$. Then \er{tseN} implies
\[
\lb{Mse}
\frac1{\n^2d}\leq\m(\wt\o_s)=\frac12\,M_s\leq{1\/2\n}\sum\limits_{\be\in\supp\gm}\lan\gm(\be),\wt\o_s\ran^2, \qqq s\in\N_d.
\]
Summing the upper estimates we have
$$
\Tr M=\sum\limits_{s=1}^d M_s\leq{1\/\n}\sum\limits_{s=1}^d\sum\limits_{\be\in\supp\gm}\lan\gm(\be),\wt\o_s\ran^2
={1\/\n}\sum\limits_{\be\in\supp\gm}\|\gm(\be)\|^2.
$$
Thus, using the inequality $M_1\leq M\leq M_d$ which is understood in the sense of quadratic forms, we get
$$
\frac2{\n^2d}\leq M_1\leq M\leq M_d\leq\Tr M\leq\frac{C_\gm}\n\,,
$$
where $C_\gm$ is defined in \er{Mtem}. Since $m=M^{-1}$, we obtain the estimate \er{Mtem}.

Let $\cI=d$. Then, due to Theorem \ref{TDImf}.\emph{i}), for some basis $a_1,\ldots,a_d$ of the lattice $\G$ the set $\{\gm(\be): \be\in\cE_*\cap\supp \gm\}$ is an orthonormal basis of $\Z^d$. In this case $C_\gm=2d$ and the lower bound on the effective mass tensor in \er{Mtem} takes the form
${\n\/2d}\leq m$. \qq \BBox

\medskip

\footnotesize
 \textbf{Acknowledgments. \lb{Sec8}}  Our study was
supported by the RSF grant  No. 18-11-00032. We would like to thank a referee for thoughtful comments that helped us to improve the manuscript.


\begin{thebibliography}{9999} \lb{Sec8}
\setlength{\itemsep}{-\parskip}\footnotesize

\bibitem[B74]{B74} Biggs, N. Algebraic graph theory,
Cambridge University Press, 1974.

\bibitem[CDS95]{CDS95} Cvetkovic, D.; Doob, M.; Sachs, H.
Spectra of graphs. Theory and applications. Johann Ambrosius Barth, Heidelberg, 1995.

\bibitem[FLP17]{FLP17} Fabila-Carrasco, J.S.; Lled\'o, F.; Post, O. Spectral gaps and discrete magnetic Laplacians, Linear Algebra Appl.,
547 (2018), no. 15, 183--216.

\bibitem[Ha02]{Ha02} Harris, P.  Carbon nano-tubes and related structure,
Cambridge, Cambridge University Press, 2002.

\bibitem[HN09]{HN09} Higuchi, Y.; Nomura, Y.
Spectral structure of the Laplacian on a covering graph.
European J. Combin. 30 (2009), no. 2, 570--585.

\bibitem[HS99a]{HS99a} Higuchi, Y.; Shirai, T. The spectrum
of magnetic Schr\"odinger operators on a graph with periodic
structure, J. Funct. Anal., 169 (1999), 456--480.

\bibitem[HS99b]{HS99b} Higuchi, Y.; Shirai, T. A remark on the
spectrum of magnetic Laplacian on a graph, the proceedings of
TGT10, Yokohama Math. J., 47 (1999), Special issue, 129--142.

\bibitem[HS04]{HS04} Higuchi, Y.; Shirai, T. Some spectral and geometric properties for infinite graphs, AMS Contemp. Math. 347 (2004), 29--56.

\bibitem[K00]{K00} Korotyaev, E. Estimates for the Hill operator. I, J.
Differential Equations 162 (2000), no. 1, 1--26.

\bibitem[K06]{K06} Korotyaev, E. Estimates for the Hill operator. II, J.
Differential Equations 223 (2006), no. 2, 229--260.

\bibitem[K08]{K08} Korotyaev, E. Effective masses for zigzag nanotubes
in magnetic fields, Lett. Math. Phys., 83 (2008), no. 1,  83--95.

\bibitem[KS14]{KS14} Korotyaev, E.; Saburova, N.
Schr\"odinger operators on periodic discrete graphs, J. Math. Anal. Appl.,  420 (2014), no. 1, 576--611.

\bibitem[KS15]{KS15} Korotyaev, E.; Saburova, N. Spectral band localization for Schr\"odinger operators on periodic  graphs, Proc. Amer. Math. Soc., 143 (2015), 3951--3967.

\bibitem[KS16]{KS16} Korotyaev, E.; Saburova, N. Effective masses
for Laplacians on periodic graphs, J. Math. Anal. Appl. 436 (2016), no. 1, 104--130.

\bibitem[KS17]{KS17}  Korotyaev, E.; Saburova, N.
Magnetic Schr\"odinger operators on periodic discrete graphs, J.
Funct. Anal., 272 (2017), 1625--1660.

\bibitem[KoS18]{KoS18}  Korotyaev, E.; Slousch, V. Asymptotics and estimates of the discrete spectrum of the Schr\"odinger operator on a discrete periodic graph, preprint, arXiv:1903.11810.

\bibitem[KSS98]{KSS98} Kotani, M.; Shirai,T.; Sunada, T.
Asymptotic behavior of the transition probability of a random walk
on an infinite graph, J. Funct. Anal. 159 (1998), no.~2, 664--689.

\bibitem[LP08]{LP08} Lled\'o, F.; Post, O. Eigenvalue bracketing for discrete and metric graphs, J. Math. Anal. Appl. 348 (2008), 806--833.

\bibitem[M91]{M91} Mohar, B. Some relations between analytic and
 geometric properties of infinite graphs, Discrete mathematics 95 (1991),
  193--219.

\bibitem[NG04]{NG04} Novoselov, K.S.; Geim, A.K. et al,
Electric field effect in atomically thin carbon films, Science 22 October,
306 (2004), no. 5696, 666--669.

\bibitem[RS78]{RS78} Reed, M.; Simon, B. Methods of modern mathematical
physics, vol.IV. Analysis of operators, Acad. Press, New York, 1978.

\bibitem[S13]{S13} Sunada, T. Topological crystallography, Surveys Tutorials  Appl. Math. Sci., vol. 6, Springer, Tokyo, 2013.

\bibitem[SS92]{SS92} Sy, P.W.; Sunada, T. Discrete Schr\"odinger
 operator on a graph, Nagoya Math. J., 125 (1992), 141--150.


\end{thebibliography}
\end{document}